\documentclass[final,leqno,letterpaper]{etna2}

\usepackage[normalem]{ulem}
\usepackage{amsmath}
\usepackage{amsfonts}
\usepackage{amssymb}
\usepackage{graphicx}
\usepackage{epstopdf,euscript}
\usepackage{url}
\usepackage{multirow}
\usepackage{array}
\usepackage{color}
\usepackage[ruled, linesnumbered,algo2e]{algorithm2e}
\usepackage{mathtools}
\usepackage{subcaption}
\usepackage{cleveref}

\usepackage{tikz}
\usetikzlibrary{quotes,angles,trees}

\usepackage[symbol]{footmisc}

\usepackage{wrapfig}
\usepackage{lscape}
\usepackage{rotating}

\def\cal{\mathcal}
\def\abs#1{\lvert#1\rvert}

\def\nm#1{\left\|#1\right\|}
\def\norm#1{\nm{#1}_{2}}

\def\C{\mathbb{C}}

\def\R{\mathbb{R}}

\def\Rn{\R^n}

\def\Rnn{\R^{n\times n}}

\def\BA{{\bf A}}  
  
  \def\CC{{\cal C}}

  \def\CK{{\cal K}}

  \def\CO{{\cal O}}
  \def\CP{{\cal P}}
  
  \def\CR{{\cal R}}

  \def\CU{{\cal U}}
  \def\CV{{\cal V}}
  \def\CW{{\cal W}}

\def\BAs{\BA{\kern-1.5pt}}

\def\CPs{\CP{\kern-0.8pt}}

\mathcode`\@="8000
{\catcode`\@=\active \gdef@{\mkern1mu}}

\def\mydate{\number\day\ {\ifcase\month \or January\or February\or
              March\or April\or May\or June\or July\or August\or
              September\or October\or November\or December\fi}
\number\year}
\def\vek#1{\mathbf{#1}}

\providecommand{\argmin}[1]{\underset{#1}{\text{{\rm argmin}}}}

\providecommand{\prn}[1]{\left(#1\right)}
\providecommand{\bigprn}[1]{\big(#1\big)}

\providecommand{\brac}[1]{\left[#1\right]}
\def\curl#1{\left\{#1\right\}}

\newcommand\restr[2]{{
  \left.\kern-\nulldelimiterspace 
  #1 
  \vphantom{\big|} 
  \right|_{#2} 
  }}

\def\Span{\rm span}

\def\be{\begin{equation}}
\def\ee{\end{equation}}
\def\bea{\begin{eqnarray}}
\def\eea{\end{eqnarray}}
\def\nn{\nonumber}
\def\mand{\mbox{\ \ \ and\ \ \ }}

\def\mwith{\mbox{\ \ \ with\ \ \ }}
\def\mwhere{\mbox{\ \ \ where\ \ \ }}

\def\msuchthat{\mbox{\ \ \ such that\ \ }}

\def\mselect{\mbox{\ \ \ select\ \ }}

\def\bbmat{\begin{bmatrix}}
\def\ebmat{\end{bmatrix}}

\def\balg{\begin{algorithm}}
\def\ealg{\end{algorithm}}
\def\balgte{\begin{algorithm2e}}
\def\ealgte{\end{algorithm2e}}
\def\bthm{\begin{theorem}}
\def\ethm{\end{theorem}}
\def\blem{\begin{lemma}}
\def\elem{\end{lemma}}
\def\bprop{\begin{proposition}}
\def\eprop{\end{proposition}}
\def\bcor{\begin{corollary}}
\def\ecor{\end{corollary}}
\def\bdefin{\begin{definition}}
\def\edefin{\end{definition}}
\def\bc{\begin{cases}}
\def\ec{\end{cases}}
\newcommand\bproof[1]{\par\addvspace{1ex} \indent\textit{Proof.}\ \ #1}
\def\eproof{\hfill\cvd\linebreak\indent}
\newtheorem{exple}{Example}
\def\bex{\begin{exple}}
\def\eex{\end{exple}}
\newtheorem{exercise}{Exercise}
\def\bexer{\begin{exercise}}
\def\eexer{\end{exercise}}
\newtheorem{conjecture}{Conjecture}
\def\bconj{\begin{conjecture}}
\def\econj{\end{conjecture}}
\newtheorem{assumption}{Assumption}
\def\bass{\begin{assumption}}
\def\eass{\end{assumption}}
\newtheorem{notation}{Notation}
\def\bnot{\begin{notation}}
\def\enot{\end{notation}}
\def\brem{\begin{remark}}
\def\erem{\end{remark}}

\def\blankfootnote#1{\let\thefootnote\relax\footnotetext{#1}}
\def\cvd{~\vbox{\hrule\hbox{%
  \vrule height1.3ex\hskip0.8ex\vrule}\hrule } }

\def\bitem{\begin{item}}
\def\eitem{\end{itemize}}
\def\benum{\begin{enumerate}}
\def\eenum{\end{enumerate}}
  
  \def\bPi{\boldsymbol\Pi}
  \def\bPhi{\boldsymbol\Phi}

\DeclareMathSymbol{\dprod}{\mathbin}{operators}{"3A}

\providecommand{\file}[1]{\texttt{\nolinkurl{#1}}}

\newcommand{\KMS}[1]{{ #1}}

\setbibdata{1}{xx}{46}{2017} 

\hypersetup{%
    pdftitle={A note on augmented unprojected Krylov subspace methods},
    pdfauthor={Kirk M. Soodhalter},
    pdfkeywords={Krylov subspaces, augmentation, recycling, discrete ill-posed problems	}
    }

\title{A note on augmented unprojected Krylov subspace methods\thanks{%
Received... Accepted... Published online on... Recommended by....
}}

\author{Kirk M. Soodhalter\footnotemark[2]}

\shorttitle{AUGMENTED UNPROJECTED KRYLOV SUBSPACES} 
\shortauthor{K.M. SOODHALTER}

\begin{document}

\maketitle

\renewcommand{\thefootnote}{\fnsymbol{footnote}}
\newcommand{\boldeta}{\boldsymbol\eta}

\footnotetext[2]{School of Mathematics, Trinity College Dublin, The University of Dublin,
College Green, Dublin 2, Ireland.
({\tt ksoodha@maths.tcd.ie}).).}

\begin{abstract}
Subspace recycling iterative methods and other subspace augmentation schemes are a successful extension to Krylov subspace methods in which
a Krylov subspace is augmented with a fixed subspace spanned by vectors deemed to be helpful in accelerating convergence or conveying knowledge
of the solution.  Recently, a survey was published, in which a framework describing the vast majority of such methods was proposed 
[Soodhalter et al, GAMM-Mitt. 2020].  In many of these methods, the Krylov subspace is one generated by the system matrix composed with a projector
that depends on the augmentation space.  However, it is not a requirement that a projected Krylov subspace be used.  There are augmentation
methods built on using Krylov subspaces generated by the original system matrix, and these methods also fit into the general framework.  

In this note, we observe that one gains implementation benefits by considering such augmentation methods with unprojected Krylov subspaces
in the general framework.  We demonstrate this 
by applying the idea to the R$^3$GMRES method proposed in [Dong et al. ETNA 2014] to obtain a simplified implementation and to connect that algorithm
to early augmentation schemes based on flexible preconditioning [Saad. SIMAX 1997].
\end{abstract}



\section{Introduction}\label{section.intro}
Augmented and recycled Krylov subspace methods have been proposed for accelerating iterative methods for solving a linear 
system (e.g., \cite{Morgan.GMRESDR.2002})
or a sequences of linear systems (see, e.g., \cite{Parks.deSturler.GCRODR.2005}) by approximating the solution to each linear system from
the sum of a Krylov subspace $\CV_{j}$ and a fixed subspace $\CU$.  
The survey \cite{dSKS.2020} details many instances of 
such methods in the literature and proposes a framework which describes their general mechanics and common mathematical structure they all share.  
In most cases, the Krylov subspace use by such a method is a \textit{projected} Krylov subspace, meaning the matrix is composed with a projector
which depends on $\CU$.  However, there are examples in the literature of augmented methods which use an unprojected Krylov subspace built
using the only the matrix, see, e.g., \cite{DGH.2014,ErhelGuyomarch:2000:1}.  Such methods necessarily also fit into the framework but are not generally
described as such.   In this note, we focus on one such method, R$^3$GMRES, 
proposed in \cite{DGH.2014}; we show how considering it as an augmented method in the
framework from \cite{dSKS.2020} allows for a simpler implementation built on well-understood algorithmic blocks from classical GMRES \cite{Saad.GMRES.1986}.
In addition, we point out that the R$^3$GMRES can be related to an older augmentation scheme built on 
flexible preconditioning \cite{Saad.Deflated-Aug-Krylov.1997}.

\section{Background}\label{section.background}
In principle, the R$^3$GMRES can be applied to any \KMS{square, discrete linear problem}, but it is proposed specifically to treat discrete ill-posed problems.
Therefore, we begin with a brief description of the ill-posed problem setting.

\KMS{%
Ill-posed problems arise often in the context of scientific applications in which one cannot directly
observe the object or quantity of interest.  
However, indirect observations or measurements can be made.	
We restrict ourselves to the linear case, whereby the unobservable quantity of interest and the measured data can be
related by a linear operator.
In this note, we consider a discretized, finite dimensional version of this problem,
\be\label{eqn.Axb}
	\vek A\vek x = \vek b \mwith \vek A\in \Rnn\mand\vek b\in\Rn.
\ee
The vector $\vek b$ represents the observed data, obtained from measurements (untainted by measurement noise), and 
$\vek x$ represents the quantity of interest, which cannot be directly observed.
The matrix $\vek A$ is generally taken to have large condition number and singular values which decrease smoothly,
with no breaks to indicate a separation between the well-posed and ill-posed parts of the discrete problem. 
As this is a finite-dimensional discretized problem, we expect perturbations in the right-hand side to produce bounded
perturbations in the reconstructed solution that are generally large enough to render the reconstructed solution useless.  
Thus, we must consider regularization methods.  In this note, we concern ourselves with some GMRES-based regularization techniques for sparse,
large-scale problems, 
but there is an extensive literature on the topics of Krylov subspace methods and hybrid methods; see, e.g., 
the surveys \cite{ChungGazzola:2021:1,GazzolaLandman:2020:1}.
}%

In the next section, we review some general mathematics behind Krylov subspace methods.  
We explain briefly GMRES before turning our attention to augmented Krylov subspace methods. 
In \Cref{section.aug}, we review the basic mechanics of augmented/recycled Krylov subspace methods, particularly in the context of
the framework proposed in \cite{dSKS.2020}. In \Cref{section.simplified-R3GMRES}, we show how the R$^3$GMRES 
method can be simplified by casting it in this framework.  Finally, we demonstrate the behavior of the new implementation with some
numerical experiments in \Cref{section.numerics}.
\bnot
In this paper, we denote by $\vek  I_{\ell}\in\R^{\ell\times \ell}$ the identity matrix acting on $\R^{\ell}$, and if the dimensions are understood from context, 
we simply write $\vek I$.  Additionally, $\underline{\vek I}_{\ell}\in\R^{(\ell+1)\times\ell}$ denotes the same identity matrix but
with an extra row of zeros appended at the bottom.  \KMS{The vector $\vek x_{0}\in\Rn$ denotes the initial approximation.}
We denote the initial error $\boldeta_{0} = \vek x - \vek x_{0}$ and the initial residual $\vek r_{0} = \vek A\boldeta_{0}  = \vek b - \vek A\vek x_{0}$.
The vector $\vek e_{i}$ denotes the $i$th canonical basis vectors whose length is defined by context.
\enot
\section{Background}\label{section.background}
In this section, we begin with a general description of Krylov subspace iterative methods, specifically the 
Generalized Minimum Residual Method (GMRES).  We then offer a brief review of augmented Krylov methods, which have been developed
both in the well-posed and ill-posed problems literature.  We observe that there has been some overlap in the 
developments in the two communities.    

\subsection{Krylov Subspace Methods}\label{section.kryl-basic}
Krylov subspace iterative methods are a well-known class of methods for the solution of linear systems 
as well as other types of problems.  For solving a linear system of the form \cref{eqn.Axb}
with $\vek A\in\Rnn$ and $\vek b\in\Rn$, one builds the Krylov subspace
\be\nn
	\CK_{j}(\vek A,\vek r_{0}) = \curl{\vek r_{0},\vek A\vek r_{0},\vek A^{2}\vek r_{0},\ldots,\vek A^{j-1}\vek r_{0}}
\ee
iteratively (at the cost of one matrix-vector product per iteration).  At iteration $j$, a correction 
$\vek t_{j}\in\CK_{j}(\vek A,\vek r_{0})$ is selected according to some constraints on the residual $\vek r_{j} = \vek b - \vek A\vek x_{j}$, where
$\vek x_{j} = \vek x_{0} + \vek t_{j}$ is the $j$th approximation.  We call $\vek t_{j}$ a \emph{correction} and the space
from which it is drawn the \emph{correction space}. 
In this paper, we focus on the Generalized Minimum
Residual Method (GMRES) \cite{Saad.GMRES.1986} 
in which we select
\be\label{eqn.pg-GMRES}
\vek t_{j}\in\CK_{j}(\vek A,\vek r_{0})\msuchthat \vek r_{j}\perp\vek A\CK_{j}(\vek A,\vek r_{0}).
\ee
Such an orthogonality condition on the residual is called a \emph{Petrov-Galerkin condition}. 
Methods with such a residual orthogonality constraint are often called \emph{residual projection methods} because
the constraint leads to a projection (oblique or orthogonal) of the residual.
This particular constraint is equivalent to solving the residual minimization problem
\be\label{eqn.gmres-full-min}
	\vek t_{j} = \argmin{\vek t\in\CK_{j}(\vek A,\vek r_{0})}\norm{\vek b - \vek A(\vek x_{0} + \vek t)},
\ee
and the residual is projected orthogonally (implicitly) to obtain the updated approximation, with
\be\label{eqn.GMRES-resid-proj}
	\vek t_{j} = \vek P_{\CK_{j}}\boldeta_{0}\mand \vek r_{j} = \prn{\vek I - \vek Q_{\CK_{j}}}\vek r_{0},
\ee
with $\vek P_{\CK_{j}}$ being the $\prn{\vek A^{\ast}\vek A}$-orthogonal projector onto $\CK_{j}\prn{\vek A,\vek r_{0}}$, and 
$ \vek Q_{\CK_{j}}$ being the orthogonal projector onto $\vek A\CK_{j}\prn{\vek A,\vek r_{0}}$.
During the iteration, one builds an orthonormal basis for $\CK_{j}(\vek A,\vek r_{0})$ one vector at-a-time 
using the Arnoldi process.  At iteration $j$, the process has generated 
\be\nn
	\vek V_{j+1}=\bbmat \vek v_{1} & \vek v_{2} & \cdots & \vek v_{j+1} \ebmat\in\R^{n\times (j+1)}\mand \underline{\vek H}_{j}\in\R^{(j+1)\times j},
\ee
where the columns of $\vek V_{j}$ form an orthonormal basis for $\CK_{j}(\vek A,\vek r_{0})$, and  
$\underline{\vek H}_{j}$ is an upper Hessenberg matrix (with zeros below the first subdiagonal) 
containing the orthogonalization coefficients.  From the 
construction of the basis, we get the Arnoldi relation
\be\label{eqn.arnoldi-relation}
	\vek A\vek V_{j} = \vek V_{j+1}\underline{\vek H}_{j} = \vek V_{j}\vek H_{j} + h_{j+1,j}\vek v_{j+1}\vek e_{j}^{T},
\ee
where $\vek H_{j}\in\R^{j\times j}$ is simply the first $j$ rows of $\underline{\vek H}_{j}$.
From \eqref{eqn.arnoldi-relation}, one can reduce the minimization \eqref{eqn.gmres-full-min} to a smaller $(j+1)\times j$
least-squares minimization problem
\be\label{eqn.gmres-ls-min}
	\vek y_{j} = \argmin{\vek y\in\R^{j}}\norm{\underline{\vek H}_{j}\vek y - \beta\vek e_{1}^{(j+1)}}\mand \vek t_{j} = \vek V_{j}\vek y_{j},
\ee
where $\beta = \norm{\vek r_{0}}$.  The standard implementation dictates that we compute the {\tt QR}-factorization
 $\underline{\vek H}_{j} = \vek Q_{j}\underline{\vek R}_{j}$ using Givens rotations, where $\vek Q_{j}\in\R^{(j+1)\times(j+1)}$
 is an orthogonal matrix,
and $\underline{\vek R}_{j}\in\R^{(j+1)\times j}$ is upper triangular.  \KMS{Via the economy {\tt QR}-factorization
of $\underline{\vek H}_{j}$}, we can 
 recast the minimization in \cref{eqn.gmres-ls-min} as the solution of an upper triangular linear system
 \be\nn
 	\vek R_{j}\vek y_{j} = \prn{\vek Q_{j}^{T}(\beta\vek e_{1})}_{1:j},
 \ee
 where $\vek R_{j}\in\R^{j\times j}$ is simply the first $j$ rows of $\underline{\vek R}_{j}$, and 
 $\prn{\cdot}_{1:j}$ denotes taking the first $j$ rows of the argument.  This can be used to develop a progressive 
 formulation of GMRES; but more importantly, it allows one to monitor the residual norm without computing the 
 GMRES approximation at each iteration.  One can show that the $j$th residual norm is simply the $(j+1)$st row of
 $\vek Q_{j}^{T}(\beta\vek e_{1})$ \cite{Saad.GMRES.1986}.

\subsection{Range-restricted Krylov subspace methods}\label{section.rrKryl}  In the context of ill-posed problems, range-restricted methods have been proposed, wherein the Krylov subspace used is $\CK_{j}(\vek A,\vek A\vek r_{0})$ rather than simply
generating with the residual $\vek r_{0}$.  The rationale in this setting is that the right-hand side (and therefore initial residual)
may be profoundly noise-polluted in such a way that reduces the effectiveness of the Krylov subspace method.  In such problems,
the matrix $\vek A$ is a discretized version of an operator that often has smoothing properties, meaning $\vek A\vek r_{0}$ is a 
smoothed version of the initial data, and using the range-restricted subspace will produce a more stable iteration.  Range-restricted
versions of GMRES \cite{RY.2005,NRS.2012,NRS-2.2012} and MINRES \cite{DMR.2014} have been proposed, with the latter being a practical realization of the MR2 method
discussed in Hanke's monograph \cite{H1995-book}.

\subsection{Augmented methods for well- and ill-posed problems}\label{section.aug-deriv}
Augmented Krylov subspace methods have been discussed in both the well- and ill-posed problems communities, though in each
with different goals in mind.  The term \emph{augmented Krylov subspace method} describes here an iterative method
in which, in addition to generating a Krylov subspace, one wishes to include vectors in the correction space 
deemed useful for \KMS{either accelerating the convergence to solution or improving the quality of the approximation delivered by
the method.  }

\KMS{%
For well-posed problems, these vectors may span a subspace
which has been determined to have strongly contributed to speed-of-convergence \cite{deSturler.GCROT.1999} or to attempt to damp the influence of 
certain parts of the spectrum of the operator \cite{Morgan.GMRESDR.2002,Parks.deSturler.GCRODR.2005}.  For ill-posed problems, this strategy
has also been shown to be effective in the case that, e.g., the noise level is rather low, as the solution may 
require many iterations \cite{dSC.2019}.

However, in the context of large-scale, discrete ill-posed problems, one may also augment with 
vectors representing known features of the image, usually those
which are highly local, such as discontinuous jumps or areas of high gradient, which an iterative method based on a Krylov subspace method may have difficulty
resolving \cite{DGH.2014,BR.2007,BR-2.2007}. Recycling-based strategies have also been shown to be effective for some such applications
\cite{dSC.2019}.  Recently, using the framework from \cite{dSKS.2020}, augmented methods were analyzed formally 
\cite{HRS.Augmented-regularization.2020} as regularization methods.
}%

\KMS{%
\subsection{Subspace augmentation via a minimization constraint}\label{section.aug}
We briefly present a general residual constraint framework through which the
methods in question can be viewed.  For a more complete view of this framework, see \cite{dSKS.2020} 
in terms of residual constraints on top of the existing work in \cite{Gaul.2014-phd,GGL.2013,Gutknecht.AugBiCG.2014,Gutknecht:2012:1}.  

In the framework, we approach augmented methods by approximating the correction over the sum of two subspaces,
$\CU$ which is \emph{fixed} and $\CV_{j}$ which is \emph{built iteratively} (i.e., it generally is some sort of Krylov subspace).
In this note, we consider the special case that we apply a residual-minimizing constraint. 
}
This technique is a straightforward generalization of the minimum residual projection
constraint \cref{eqn.pg-GMRES}, i.e., we require 
\be\label{eqn.aug-min-constraint}
	\vek b - \vek A\prn{\vek x_{0} + \vek s_{j} + \vek t_{j}}\perp\vek A\prn{\CU +\CV_{j}}.
\ee
This residual constraint underpins (either implicitly or explicitly) many augmented GMRES-type methods.  
Associated to this constraint \KMS{are, respectively, the
$\vek A^{T}\vek A$-orthogonal projector onto $\CU$ and the orthogonal projector onto $\vek A\,\CU$}
\be\nn
	\bPi = \vek U\prn{\vek U^{T}\vek A^{T}\vek A\vek U}^{-1}\vek U^{T}\vek A^{T}\vek A\mand \bPhi = \vek A\vek U\prn{\vek U^{T}\vek A^{T}\vek A\vek U}^{-1}\vek U^{T}\vek A^{T}.
\ee
If a method minimizes the residual over a sum of subspaces, it \emph{necessarily} fits
into the augmentation framework, regardless of how $\CV_{j}$ is generated.  
\KMS{It is pointed out in \cite{dSKS.2020} that regardless of the choice of $\CV_{j}$, this residual minimization over the sum
of subspaces can be reduced and reformulated as selecting $\vek t\approx\vek t_{j} \in\CV_{j}$ to minimize the residual of the projected problem
\begin{align*}
	(\vek I-\bPhi)\vek A\vek t = (\vek I-\bPhi)\vek b
\end{align*}
and setting $\vek x_{j} = \vek x_{0} + \bPi\boldsymbol\eta_{0} + (\vek I - \bPi)\vek t_{j}$, where we note that the action of $\bPi$ on $\boldsymbol\eta_{0}$ can be
computed efficiently without knowing $\boldsymbol\eta_{0}$.}
We show that this can lead to simplified
implementations of such methods, particularly as it relates to methods which augment \emph{unprojected} Krylov subspace methods.

For methods such as GMRES-DR \cite{Morgan.GMRESDR.2002} and GCRO-type methods, e.g.,\cite{deSturler.GCROT.1999,deSturler.GCRO.1996,Parks.deSturler.GCRODR.2005}, 
the iteratively generated Krylov subspace matches with the projected subproblem \cref{eqn.projAxb}
with $\CV_{j} = \CK_{j}\prn{\prn{\vek I - \bPhi}\vek A,\prn{\vek I - \bPhi}\vek r_{0}}$.  
Augmented methods based on range-restricted GMRES, 
e.g., \cite{BR.2007}, use $\CV_{j} = \CK_{j}\prn{\prn{\vek I - \bPhi}\vek A,\prn{\vek I - \bPhi}\vek A\vek r_{0}}$.  
With GCRO-based augmented (range-restricted) GMRES, one can implement either one small minimization problem over the
augmented subspace or by directly using the above framework to approximate the solution of \cref{eqn.projAxb} by a GMRES minimization followed by a projection,
as described above.  Let  
\be\nn
\CV_{j} = \CK_{j}\prn{\prn{\vek I - \bPhi}\vek A,\prn{\vek I - \bPhi}\vek w_{0}}\mwhere \vek w_{0} \in\curl{ \vek r_{0},\,\vek A\vek r_{0}}.
\ee
Let $\vek V_{j}$ be generated by the Arnoldi process so that we have
\be\label{eqn.projAxb}
	\prn{\vek I - \bPhi}\vek A\vek V_{j} = \vek V_{j+1}\underline{\vek H}_{j}.
\ee
\bass\label{assumpt.Corth}
The matrix $\vek U$ is scaled so that $\vek C = \vek A\vek U$ has orthonormal columns, i.e., $\vek C^{T}\vek C=\vek I_{k}$ and $\bPhi=\vek C\vek C^{T}$.
This is not mathematically necessary, but it allows for various algorithmic simplifications.
\eass

In \cite{Parks.deSturler.GCRODR.2005}, the authors approach recycling by deriving a modified Arnoldi relation
\small\be\label{eqn.mod-arnoldi}
	\vek A\bbmat \vek U & \vek V_{j} \ebmat = \bbmat \vek C & \vek V_{j+1} \ebmat \underline{\vek G}_{j}\mwhere  \underline{\vek G}_{j}=\bbmat \vek I_{k} & \vek B_{j}  \\ & \underline{\vek H}_{j} \ebmat\mand \vek B_{j} = \vek C^{T}\vek A\vek V_{j}.
\ee\normalsize
From this, one can satisfy \cref{eqn.aug-min-constraint} by solving the small least squares problem
\be\label{eqn.aug-GMRES-ls}
	\prn{\vek z_{j}, \vek y_{j}}=\argmin{ \vek u\in\R^{k}\atop \vek v \in\R^{j}}\norm{ \bbmat \vek C & \vek V_{j+1} \ebmat^{T}\vek r_{0} - \underline{\vek G}_{j}\bbmat \vek z\\\vek y  \ebmat }.
\ee
This is in actuality not necessary for implementing the method, since it can be decoupled to solve a GMRES small least-squares 
problem for $\vek y_{j}$
which then enables the solution of $\vek z_{j}$ by back substitution; i.e., 
\begin{equation}\label{eqn.gcro-yz-calcs}
	\vek y_{j} = \argmin{\vek y\in\C^{j}}\nm{\beta\vek e_{1} - \underline{\vek H}_{j}\vek y}\mand \vek z_{j} = \vek C^{\ast}\vek r_{0}-\vek B_{j}\vek y_{j}.
\end{equation}
\KMS{However, \cref{eqn.aug-GMRES-ls} is useful as a comparison to the coupled minimization
in the proposed implementation of augmented unprojected (range-restricted) GMRES, which we discuss below.}

\subsection{Augmenting unprojected Krylov subspaces}
In both the well- and ill-posed problems community, augmented methods have been 
proposed wherein an unprojected Krylov subspace is used, in
\cite{Saad.Deflated-Aug-Krylov.1997} $\CV_{j} = \CK_{j}\prn{\vek A,\vek r_{0}}$ and in \cite{DGH.2014} $\CV_{j} = \CK_{j}\prn{\vek A,\vek A\vek r_{0}}$. We discuss briefly some implementation details of these methods \emph{which are relevant to the present note}, but one should read the cited papers and
references therein for complete implementation details.  It has also been observed that under certain circumstances in which there are strict
constraints on the amount of computations one can perform per iteration, an unprojected augmented method may be preferred (or indeed be the only option);
see \cite{RamlauStadler:2020:1}, which builds on \cite{ErhelGuyomarch:2000:1}.

\subsubsection{Flexible GMRES-based augmentation for well-posed problems}
In \cite{Saad.Deflated-Aug-Krylov.1997}, Saad proposes augmenting an already constructed Krylov subspace 
$\CK_{j}(\vek A,\vek r_{0})$ with a subspace $\CW=\Span\curl{\vek w_{1},\vek w_{2},\ldots,\vek w_{k}}$ by 
treating the basis vectors of $\CW$ as those resulting from the action of successive implicit flexible preconditioners. 
The augmentation process is embedded in an iteration of flexible GMRES.  This minimum residual method can be
described in the language of the framework by identifying that the correction space in this setting is
$\underbrace{\CK_{j}(\vek A,\vek r_{0})}_{\CV_{j}} + \underbrace{\CW}_{\CU}$, and the constraint space is 
$\underbrace{\vek A\CK_{j}(\vek A,\vek r_{0})}_{\widetilde{\CV}_{j}} + \underbrace{\vek A\CW}_{\widetilde{\CU}}$, where
the flexible Arnoldi process produces an orthonormal basis for the constraint space. An outline of this method
is shown in Algorithm \ref{alg.fgmres-aug}.

\RestyleAlgo{boxruled}
\balgte[h]
\caption{One cycle of Flexible GMRES-based augmentation from \cite{Saad.Deflated-Aug-Krylov.1997}\label{alg.fgmres-aug}}
\SetKwInOut{Input}{Input}\SetKwInOut{Output}{Output}
\Input{$\vek A\in\Rnn$, $\vek b,\vek x_{0}\in\Rn$, $\vek W\in\R^{n\times k}$, $\vek m>0$}
$\vek r_{0} = \vek b - \vek A\vek x_{0}$; $\beta = \nm{\vek r_{0}}$; $\vek v_{1} = \vek r_{0}/\beta$\;
\For{$i=1,2,\ldots, \KMS{m+k}$}{
	\uIf{$i< m$}{
		$\vek v_{i+1} = \vek A\vek v_{i}$\;
	}\Else{
		$\vek v_{i+1} = \vek A\vek w_{i-m+1}$\;
	}
	\For{$j=1,2,\ldots i$}{
		$h_{ji} = \vek v_{j}^{\ast}\vek v_{i+1}$\;
		$\vek v_{i+1}\leftarrow \vek v_{i+1} - h_{ji}\vek v_{j}$\;
	}
	$h_{i+1,i} = \nm{\vek v_{i+1}}$\;
	$\vek v_{i+1}\leftarrow \vek v_{i+1}/h_{i+1,i}$\;
	$\vek y = \argmin{\vek y\in\R^{i}}\nm{\beta\vek e_{1} - \underline{\vek H}_{i}\vek y}$\;
}
\KMS{$\vek x = \vek x_{0} + \vek V_{m}\vek y{\tt (1:m)} + \vek W\vek y{\tt(m+1:m+k)}$}\;
\ealgte
\subsubsection{Augmentation of (range-restricted) methods for solving ill-posed problems}
In the context of solving discrete ill-posed problems using augmented iterative techniques, it has been asserted 
in \cite{DGH.2014} that it may be preferable to 
employ augmentation techniques with an unprojected Krylov subspace. 
This is in part motivated by the use of projected Krylov subspaces in
\cite{BR-2.2007}.  
The authors argue that the subspace $\CU$ should contain (approximations of) known features of the image.  
However, if these features are poor approximations of 
image features (e.g., a misplaced discontinuity), it is asserted that the use of a projected Krylov subspace can cause the iteration to semi-converge to a poor quality solution.  Conversely, for 
solving a well-posed problem, the iteration would eventually recover and converge.   The authors suggest using 
$\CV_{j} = \CK_{j}\prn{\vek A,\vek w_{0}}$, with $\vek w_{0}\in\curl{\vek r_{0},\vek A\vek r_{0}}$ preferring 
$\vek w_{0} = \vek A\vek r_{0}$ (i.e., a range-restricted method) as it tends to yield superior performance for their experiments
\cite{DGH.2014}. 

\brem
For the case $\vek w_{0} = \vek r_{0}$, we observe that
the method shown in \cite{DGH.2014} is \KMS{mathematically} equivalent to the augmentation in a flexible preconditioning framework proposed by 
Saad \cite{Saad.Deflated-Aug-Krylov.1997}, 
an equivalence noted in \cite{dSKS.2020}. 
\erem

\section{\color{blue}Framework perspective allows for a simplified unprojected augmented GMRES}\label{section.simplified-R3GMRES}

 Again, for implementation purposes, we invoke \Cref{assumpt.Corth}, i.e., that
 $\vek C = \vek A\vek U$ has orthonormal columns.  At each iteration of the Arnoldi process for $\CK_{j}(\vek A,\vek w_{0})$, 	
the method proposed in \cite{DGH.2014} requires an orthonormal basis for the columns of
$\vek A\bbmat \vek V_{j} & \vek U\ebmat$.  Unlike with GCRO-type methods, this does not come for free since the the Krylov subspace
is unprojected.  We show in the following subsection that approaching this method from the framework point-of-view allows us to avoid
the algorithmic complication of this orthogonalization. The framework enables us to solve for least-squares approximate solutions of the
projected problem \cref{eqn.projAxb} over the unprojected Krylov subspace $\CK_{j}\prn{\vek A,\vek w_{0}}$ and then subsequently obtain
an additional correction over $\CU$ to obtain the full approximation without additional orthogonalization complications.  Furthermore, this new
formulation allows for the estimation of the residual norm, meaning that similar to an efficient implementation of GMRES, neither the full 
approximation nor the residual need to be computed until possible convergence has been detected.

We derive a simplified version of  R$^3$GMRES in \cite{DGH.2014}. 
We begin our derivation similar to that in \cite{DGH.2014} by assuming one must progressively orthogonalize $\vek C$ against the 
Arnoldi vectors, but through our derivation we show this is actually not necessary.

Let    
\be\nn
	\widehat{\vek C}_{1} = \vek C - \vek v_{1}\prn{\vek v_{1}^{T}\vek C} \mand \widehat{\vek C}_{1} = \vek C_{1}\vek F_{1} \ \mbox{(skinny {\tt QR}-factorization)}.
\ee
At each iteration $i$, this orthogonalization must be updated after $\vek v_{i+1}$ has been generated, and this can be performed recursively
\be\nn
	\widehat{\vek C}_{i+1} = \vek C_{i} - \vek v_{i+1}\prn{\vek v_{i+1}^{T}\vek C}_{i} \mand \widehat{\vek C}_{i+1} = \vek C_{i+1}\widehat{\vek F}_{i+1} \ \mbox{(skinny {\tt QR}-factorization)}.
\ee
From this, one gets the new modified Arnoldi factorization
\be\label{eqn.rrr-mod-arnoldi}
	\vek A\bbmat \vek V_{j} & \vek U\ebmat = \bbmat \vek V_{j+1} & \vek C_{j} \ebmat\widehat{\vek G}_{j} \mwith \widehat{\vek G}_{j} =  \bbmat \underline{\vek H}_{j}& \vek D_{j} \\ & \vek F_{j}\ebmat,
\ee
where $\vek D_{j} = \vek V_{j+1}^{T}\vek C$ and $\vek F_{j} = \vek C_{j}^{T}\vek C$.  One observes that $\vek D_{j}$ can be constructed iteratively, as
\be\nn
	\vek D_{j} = \vek V_{j+1}^{T}\vek C = \bbmat \vek V_{j}^{T}\vek C \\ \vek v_{j+1}^{T}\vek C\ebmat = \bbmat \vek D_{j-1} \\ \vek d_{j}\ebmat,
\ee
where $\vek d_{j} = \vek v_{j+1}^{T}\vek C$.  

As with GCRO-based methods, the minimization constraint \cref{eqn.aug-min-constraint} reduces to a small least-squares problem similar to \cref{eqn.aug-GMRES-ls}, namely
\be\label{eqn.rrrGMRES-small-ls}
	\prn{\vek z_{j}, \vek y_{j}}=\argmin{ \vek z\in\R^{k}\atop \vek y \in\R^{j}}\norm{ \bbmat \vek V_{j+1} & \vek C_{j} \ebmat^{T}\vek r_{0} - \widehat{\vek G}_{j}\bbmat \vek y\\\vek z  \ebmat }.
\ee
This is the minimization that is then explicitly solved in \cite{DGH.2014}.  However, just like the GCRO-based methods, this minimization over the sum of subspaces
can be rewritten as the approximation of the solution of a projected subproblem \cref{eqn.projAxb} over $\CV_{j}$ whose solution is then projected onto $\CU$
to get $\vek s_{j}$.  The difference here is that $\CV_{j}$ is the Krylov subspace associated to the unprojected problem;
i.e., $\CV_{j} = \CK_{j}(\vek A,\vek w_{0})$.

The method proposed in \cite{DGH.2014} is a residual minimization over the sum of two spaces; thus it must fit into the framework 
introduced in \Cref{section.aug}.
Our task is to understand how this residual minimization over the sum of two spaces can be rewritten as a minimization for a projected subproblem, just as we have
described for GCRO-based methods.  This brings us to the main result,

\bthm\label{thm.minimizing-problem}
	Let $\CV_{j}=\CK_{j}(\vek A,\vek w_{0})$ where $\vek w_{0}\in\curl{\vek r_{0},\vek A\vek r_{0}}$.  Minimizing the residual over
	the sum of spaces $\CU + \CV_{j}$ as described in both \cite{Saad.Deflated-Aug-Krylov.1997} and \cite{DGH.2014} is equivalent
	to computing $\vek y_{j}$ satisfying
	\begin{equation}\label{aug-noproj-gmres-normeq}
		\prn{\underline{\vek H}_{j}^{T}\underline{\vek H}_{j} 
		-
		\underline{\vek H}_{j}^{T}\vek D_{j}\vek D_{j}^{T}\underline{\vek H}_{j}}\vek y_{j}
		= 
		\underline{\vek H}_{j}^{T}\prn{\vek V_{j+1}^{T}\prn{\vek I - \bPhi}\vek r_{0}}
	\end{equation}
	and $\vek z_{j}=\vek C^{T}\vek r_{0}- \vek D_{j}^{T} \underline{\vek H}_{j}\vek y_{j}$.  Furthermore, this is equivalent to finding $\vek t_{j} = \vek V_{j}\vek y_{j}\in\CV_{j}$ 
	which satisfies a least-squares minimization applied to the projected
	subproblem \cref{eqn.projAxb}, namely
	\begin{equation}\label{eqn.least-squares-thm}
		\mselect\ \vek t_{j}\in\CV_{j}\msuchthat \nm{\prn{\vek I-\bPhi}\bigprn{\vek b - \vek A\prn{\vek x_{0} + \vek t_{j}}}}\ \mbox{is minimized}.
	\end{equation}
\ethm
\bproof
One can take a couple of different approaches to see how one solves the projected subproblem.  Here we follow the approach in \cite{Parks.Soodhalter.Szyld.16}, 
wherein we form the normal equations of \cref{eqn.rrrGMRES-small-ls}
\bea
	\bbmat  \vek F_{j}^{T} & \vek D_{j}^{T}\\ & \underline{\vek H}_{j}^{T} \ebmat \bbmat \vek F_{j} & \\ \vek D_{j} & \underline{\vek H}_{j} \ebmat \bbmat   \vek z_{j} \\ \vek y_{j}\ebmat &=& \bbmat  \vek F_{j}^{T} & \vek D_{j}^{T}\\ & \underline{\vek H}_{j}^{T} \ebmat \bbmat \vek C^{T}\vek r_{0}\\ \vek V_{j+1}^{T}\vek r_{0} \ebmat\iff\nn\\
	\bbmat\vek F_{j}^{T}\vek F_{j}  + \vek D_{j}^{T}\vek D_{j} &  \vek D_{j}^{T} \underline{\vek H}_{j} \\ \underline{\vek H}_{j}^{T} \vek D_{j} &\underline{\vek H}_{j}^{T}\underline{\vek H}_{j} \ebmat \bbmat   \vek z_{j} \\ \vek y_{j}\ebmat &=& \bbmat \prn{\vek D_{j}^{T}\vek V_{j+1}^{T} + \vek F_{j}^{T}\vek C_{j}^{T}}\vek r_{0}\\ \underline{\vek H}_{j}^{T}\vek V_{j+1}^{T}\vek r_{0} \ebmat\label{eqn.rrrGMRES-normalEqn}
\eea
A block LU-factorization of the system matrix allows us to eliminate $\vek z_{j}$ from the second equation, yielding the equations
\bea
\prn{\vek F_{j}^{T}\vek F_{j}  + \vek D_{j}^{T}\vek D_{j} }\vek z_{j} +  \vek D_{j}^{T} \underline{\vek H}_{j}\vek y_{j}&=& \prn{\vek D_{j}^{T}\vek V_{j+1}^{T} + \vek F_{j}^{T}\vek C_{j}^{T}}\vek r_{0}\mand\nn\\
\prn{\underline{\vek H}_{j}^{T}\underline{\vek H}_{j} -\underline{\vek H}_{j}^{T}\vek D_{j}\prn{\vek D_{j}^{T}\vek D_{j} + \vek F_{j}^{T}\vek F_{j}}^{-1}\vek D_{j}^{T}\underline{\vek H}_{j}}\vek y_{j}&=& \underline{\vek H}_{j}^{T}\vek V_{j+1}^{T}\vek r_{0}\nn\\ - \underline{\vek H}_{j}^{T}\vek D_{j}\left(\vek D_{j}^{T}\vek D_{j} \right.&+&\left. \vek F_{j}^{T}\vek F_{j}\right)^{-1}\prn{\vek D_{j}^{T}\vek V_{j+1}^{T} + \vek F_{j}^{T}\vek C_{j}^{T}}\vek r_{0}.\nn
\eea
Observe now that if we substitute the definitions of $\vek D_{j}$ and $\vek F_{j}$ into the latter equations, we get that
\be\nn
	\vek D_{j}^{T}\vek D_{j} + \vek F_{j}^{T}\vek F_{j} = \vek C^{T}\prn{\vek V_{j+1}\vek V_{j+1}^{T} + \vek C_{j}\vek C_{j}^{T}}\vek C.
\ee
By design, we have that $\CR\prn{\vek V_{j+1}\vek V_{j+1}^{T}\vek C}\oplus \CR\prn{\vek C_{j}} = \CC$ which implies that 
\be\label{eqn.proj-is-ident}
\prn{\vek V_{j+1}\vek V_{j+1}^{T} + \vek C_{j}\vek C_{j}^{T}}\vek C = \vek C, 
\ee
and by \cref{assumpt.Corth}, we get 
\begin{align}
\vek z_{j}  &=& \prn{\vek D_{j}^{T}\vek V_{j+1}^{T} + \vek F_{j}^{T}\vek C_{j}^{T}}\vek r_{0}- \vek D_{j}^{T} \underline{\vek H}_{j}\vek y_{j}\mand\label{eqn.block-LU-decoupled}\\
\prn{\underline{\vek H}_{j}^{T}\underline{\vek H}_{j} -\underline{\vek H}_{j}^{T}\vek D_{j}\vek D_{j}^{T}\underline{\vek H}_{j}}\vek y_{j}&=& \underline{\vek H}_{j}^{T}\vek V_{j+1}^{T}\vek r_{0} - \underline{\vek H}_{j}^{T}\vek D_{j}\prn{\vek D_{j}^{T}\vek V_{j+1}^{T} + \vek F_{j}^{T}\vek C_{j}^{T}}\vek r_{0}.\nn
\end{align}
We finish by showing that
the second set of equations are the normal equations for the least-squares problem \cref{eqn.least-squares-thm} from the
statement of the theorem.  One sees this by noting that 
\be\nn
	\prn{\underline{\vek H}_{j}^{T}\underline{\vek H}_{j} -\underline{\vek H}_{j}^{T}\vek D_{j}\vek D_{j}^{T}\underline{\vek H}_{j}} = \underline{\vek H}_{j}^{T}\prn{\vek I - \vek D_{j}\vek D_{j}^{T}}\underline{\vek H}_{j} = \underline{\vek H}_{j}^{T}\vek V_{j+1}^{T}\prn{\vek I - \vek C\vek C^{T}}\vek V_{j+1}\underline{\vek H}_{j},
\ee
and that $\vek C\vek C^{T} = \bPhi$.
For the right-hand side, one observes from the definition of $\vek D_{j}$ that 
\be\nn
	\underline{\vek H}_{j}^{T}\vek D_{j}\prn{\vek D_{j}^{T}\vek V_{j+1}^{T} + \vek F_{j}^{T}\vek C_{j}^{T}}= \vek V_{j}^{T}\vek A^{T}\vek C\vek C^{T}\prn{\vek V_{j+1}\vek V_{j+1}^{T} + \vek C_{j}\vek C_{j}^{T}}.
\ee
As we have seen in \cref{eqn.proj-is-ident},
\be\nn
	\vek C^{T}\prn{\vek V_{j+1}\vek V_{j+1}^{T} + \vek C_{j}\vek C_{j}^{T}} = \brac{\prn{\vek V_{j+1}\vek V_{j+1}^{T} + \vek C_{j}\vek C_{j}^{T}}\vek C}^{T} = \vek C^{T}
\ee
which means the right-hand side of the second equation of \eqref{eqn.block-LU-decoupled} can be simplified \linebreak as
$\underline{\vek H}_{j}^{T}\vek V_{j+1}^{T}\prn{\vek I - \bPhi}\vek r_{0}$. Thus we can rewrite \cref{eqn.block-LU-decoupled} yielding
\begin{align}
\vek z_{j}  &= \vek C^{T}\vek r_{0}- \vek D_{j}^{T} \underline{\vek H}_{j}\vek y_{j}\mand\label{eqn.block-LU-updated}\\
\underline{\vek H}_{j}^{T}\vek V_{j+1}^{T}\prn{\vek I - \bPhi}\vek V_{j+1}\underline{\vek H}_{j}\vek y_{j}&= \underline{\vek H}_{j}^{T}\vek V_{j+1}^{T}\prn{\vek I - \bPhi}\vek r_{0}.\nn
\end{align}
Observing that the idempotency of projectors means $\prn{\vek I - \bPhi}=\prn{\vek I - \bPhi}^{2}$ completes the proof, since
$\bPhi$ being an orthogonal projector means it is symmetric.
\eproof

\brem
	We note that this result indicates the matrices $\vek C_{j}$ and $\vek F_{j}$ are not needed to implement R$^{3}$GMRES, greatly
	simplifying the method, as it is no longer required to progressively orthogonalize $\vek C$ with respect to to the Arnoldi vectors.
\erem

The final step in developing an efficient implementation of R$^{3}$GMRES is to rewrite and simplify
\cref{aug-noproj-gmres-normeq} using the standard Givens-rotation-based progressive {\tt QR}-factorization of 
$\underline{\vek H}_{j}$ which then enables the estimation of the residual norm without needing to 
compute the solution to \cref{aug-noproj-gmres-normeq} at each iteration.  Unlike GMRES or the GCRO-
variants of augmented methods, an exact residual norm is not available without computing 
the residual itself, which we would like to avoid.
\bthm\label{thm.compute-coeffs}
	Let $\underline{\vek H}_{j} = \vek Q_{j}\underline{\vek R}_{j}$ be the {\tt QR}-factorization obtained obtained progressively using 
	Givens rotations.  Then we can represent the coefficient vectors $\vek y_{j}$ as the solution of the linear system
	\be\label{eqn.normaleqs-subsystem}
		\prn{\vek I - \vek M_{j}\vek M_{j}^{T}}\vek R_{j}\vek y_{j} = \curl{\vek Q_{j}^{T}\prn{\vek V_{j+1}^{T}\prn{\vek I - \bPhi}\vek r_{0}}}_{1:j}
	\ee
	where $\vek M_{j} = \curl{\vek Q_{j}^{T}\vek D_{j}}_{1:j}\in\R^{j\times k}$.  
	Furthermore, the residual norm satisfies
	\begin{align}
		\nm{\prn{\vek I-\vek C\vek C^{T}}\bigprn{\vek b - \vek A\prn{\vek x_{0} + \vek V_{j}\vek y_{j}}}}^{2}&\leq \abs{\vek e^{T}_{j+1}\vek Q_{j}^{T}\vek V_{j+1}^{T}\vek r_{0}}^{2}
		\label{eqn.resid-norm-estimate}\\
		&\ \  
		+ \nm{\prn{\vek I - \vek V_{j+1}\vek V_{j+1}^{T}}\vek r_{0}}^{2},\nn
	\end{align}
	a bound which can be updated progressively.
\ethm

\noindent
\bproof
	Consider the {\tt QR}-factorization $\underline{\vek H}_{j}=\vek Q_{j}\underline{\vek R}_{j}$, obtained progressively using Givens rotations.  As has been observed in the derivation of GMRES in \cite{Saad.GMRES.1986}, we can write 
	\be\nn
		\underline{\vek H}_{j}^{T}\underline{\vek H}_{j} = \underline{\vek R}_{j}^{T}\underline{\vek R}_{j} = \vek R_{j}^{T}\vek R_{j}.
	\ee
	With this, we can rewrite \cref{aug-noproj-gmres-normeq} as
	\be\label{eqn.noproj-ne-workings1}
		\prn{\vek R_{j}^{T}\vek R_{j} - \underline{\vek R}_{j}^{T}\vek Q_{j}^{T}\vek D_{j}\vek D_{j}^{T}\vek Q_{j}\underline{\vek R}_{j}}\vek y_{j} = \underline{\vek R}_{j}^{T}\vek Q_{j}^{T}\prn{\vek V_{j+1}^{T}\prn{\vek I - \bPhi}\vek r_{0}}.
	\ee	
	We assume that Arnoldi process has not broken down and thus $\vek R_{j}$ is nonsingular.  Thus, we can multiply the 
	\cref{eqn.noproj-ne-workings1} by $\vek R_{j}^{-T}$, yielding
	\be\label{eqn.noproj-ne-workings2}
		\prn{\vek R_{j} - \underline{\vek I}_{j}^{T}\vek Q_{j}^{T}\vek D_{j}\vek D_{j}^{T}\vek Q_{j}\underline{\vek R}_{j}}\vek y_{j} = \underline{\vek I}_{j}^{T}\prn{\vek Q_{j}^{T}\prn{\vek V_{j+1}^{T}\prn{\vek I - \bPhi}\vek r_{0}}}.
	\ee
	Let $\vek M_{j} = \underline{\vek I}_{j}^{T}\vek Q_{j}^{T}\vek D_{j}\in\R^{j\times k}$.  We can 
	simplify \cref{eqn.noproj-ne-workings2} by substituting in $\vek M_{j}$, which yields
	\be\label{eqn.simplified-projected-system}
		\prn{\vek I - \vek M_{j}\vek M_{j}^{T}}\vek R_{j}\vek y_{j} = \curl{\vek Q_{j}^{T}\prn{\vek V_{j+1}^{T}\prn{\vek I - \bPhi}\vek r_{0}}}_{1:j}.
	\ee
	If the rank-$k$ outer product 
	$\vek M_{j}\vek M_{j}^{T}$ does not have any unit eigenvalues then $\prn{\vek I - \vek M_{j}\vek M_{j}^{T}}$ is invertible.  We note that this is indeed the
	case since \cref{eqn.simplified-projected-system} is derived from normal equations that have a unique solution \KMS{in this case}.
	
	Recall from the proof of \Cref{thm.minimizing-problem} that 
	the solution to this linear system $\vek y_{j}$ is the minimizer of $\norm{(\vek I - \vek C\vek C^{T})\prn{\vek V_{j+1}\underline{\vek H}_{j}\vek y - \vek r_{0}}}$.  
	As $\vek I-\vek C\vek C^{T}$ is an orthogonal projection, its action either has no effect on the vector
	norm or it reduces the length.  Thus we can estimate from above by disregarding the projector.	
	Furthermore, this analysis should include the case that the Krylov 
	subspace is range-restricted; thus $\vek r_{0}$ may not be in $\CR\prn{\vek V_{j+1}}$.  	
	As it has been pointed out
	in (e.g., \cite{NeumanReichelSadok:2012}) it suffices in this case to split the residual into $\vek V_{j+1}\vek V_{j+1}^{T}\vek r_{0}$ and the part
	in the orthogonal complement and to consider the minimization only on the part in the Krylov subspace.
	Thus we can write
	\begin{align*}
		\nm{(\vek I - \vek C\vek C^{T})\prn{\vek V_{j+1}\underline{\vek H}_{j}\vek y - \vek r_{0}}}^{2}  
		&
		\leq
		\nm{\prn{\vek V_{j+1}\underline{\vek H}_{j}\vek y - \vek r_{0}}}^{2}  
		\\
		&
		=
		\nm{\prn{\vek V_{j+1}\underline{\vek H}_{j}\vek y - \vek V_{j+1}\vek V_{j+1}^{T}\vek r_{0}}
		- \prn{\vek I -\vek V_{j+1}\vek V_{j+1}^{T} }\vek r_{0}}^{2}  
		\\
		&
		=
		\nm{\prn{\vek V_{j+1}\underline{\vek H}_{j}\vek y - \vek V_{j+1}\vek V_{j+1}^{T}\vek r_{0}}}^2
		+
		\nm{\prn{\vek I -\vek V_{j+1}\vek V_{j+1}^{T} }\vek r_{0}}^{2}  
		\\
		&
		=
		\nm{\prn{\vek Q_{j}\underline{\vek R}_{j}\vek y - \vek V_{j+1}^{T}\vek r_{0}}}^2
		+
		\nm{\prn{\vek I -\vek V_{j+1}\vek V_{j+1}^{T} }\vek r_{0}}^{2}  
		\\
		&
		=
		\nm{\prn{\underline{\vek R}_{j}\vek y - \vek Q_{j}^{T}\vek V_{j+1}^{T}\vek r_{0}}}^2
		+
		\nm{\prn{\vek I -\vek V_{j+1}\vek V_{j+1}^{T} }\vek r_{0}}^{2}  
	\end{align*}
	The result follows from the same logic used to derive the GMRES residual monitoring strategy shown in, e.g., 
	\cite{Saad.GMRES.1986}.
\eproof

\KMS{%
We note that if $\vek w_{0}=\vek r_{0}$ then 
\begin{equation}\nonumber
\nm{\prn{\vek I - \vek V_{j+1}\vek V_{j+1}^{T}}\vek r_{0}}=0.
\end{equation}
If $\vek w_{0}=\vek A\vek r_{0}$ one can progressively update and monitor this quantity by projecting $\vek r_{0}$ away
from $\CK_{j}(\vek A,\vek A\vek r_{0})$. 
}%

We observe that the estimate of the residual norm is simply the residual norm one would 
obtain from applying non-augmented (range-restricted) GMRES to the problem.  Thus, depending on the effectiveness
of the augmentation, it will likely overestimate the true residual norm.  However, the residual norm estimate \cref{eqn.resid-norm-estimate} can be
used in early iterations to avoid computing the \KMS{solution and the} residual until the estimate indicates convergence may be imminent.  
The strategy we advocate here is to use the ratio $\nm{\vek r_{0}}/\nm{\prn{\vek I-\vek C\vek C^{T}}\vek r_{0}}$
as a scaling factor between the estimate of the norm and the actual norm. This scaling factor can be updated any time the code
does an explicit residual computation, in the case we find that the estimate has falsely predicted convergence.

The matrix $\vek M_{j}$ can be constructed progressively using Givens rotations.  We initialize $\vek m_{1} = \vek d_{1}$,
reminding the reader that we are indexing the rows of $\vek M_{j}$.
At iteration $j$, we set $\vek m_{j+1}=\vek d_{j+1}$ and use the $j$th Givens rotations to make the update 
\begin{align*}
	\begin{bmatrix}
		\vek m_{j}
		\\
		\vek m_{j+1}
	\end{bmatrix}
	\leftarrow
	\begin{bmatrix}
		c_{j}
		&
		s_{j}
		\\
		-s_{j}
		&
		c_{j}
	\end{bmatrix}
	\begin{bmatrix}
		\vek m_{j}
		\\
		\vek m_{j+1}
	\end{bmatrix}
\end{align*}

We bring all this together to present a simplified implementation of R$^3$GMRES in Algorithm \ref{alg.r3gmres}.
Note that following from the strategy advocated by de Sturler \cite{EdS-personal.2020}, we compute the {\tt QR}-factorization 
$\vek A\widehat{\vek U} = \vek C\vek F$, but we do not update $\vek U=\widehat{\vek U}\vek F^{-1}$.  For $\vek z\in\R^{k}$, it is generally cheaper when
expanding $\vek U\vek z$ to calculate $\widehat{\vek U}(\vek F^{-1}\vek z)$.  This is what we do in our implementation.

\balgte
\caption{A simplified R$^3$GMRES implementation (with range restriction)\label{alg.r3gmres}}
\SetKwInOut{Input}{Input}\SetKwInOut{Output}{Output}
\Input{$\vek A\in\Rnn$, $\vek x_{0},\vek b\in\Rn$, $\vek U\in\R^{n\times k}$, $\varepsilon_{tol} > 0$}
$[\vek C,\vek F] = {\tt QR(\vek A\vek U)}$\\
$\vek r_{0} = \vek b - \vek A\vek x_{0}$; $\vek w_{0} = \vek A\vek r_{0}$\\
$\gamma \leftarrow \nm{\vek r_{0} - \vek C\vek C^{T}\vek r_{0}}/\nm{\vek r}_{0}$\\
$\vek v_{1} \leftarrow \vek w_{0}/\norm{\vek w_{0}}$\\
$\vek d_{1}\leftarrow \vek v_{1}^{T}\vek C$\\
$\vek m_{1} \leftarrow \vek d_{1}$\\ 
$\vek s_{1} = \vek U \prn{ \vek F^{-1} \prn{\vek C^{T} \vek r_{0}} }$\\
\For{$i=1,2,\ldots, j$}{
	$ \vek v_{i+1} \leftarrow \vek A\vek v_{i}$\\
	\For{$m=1,2,\ldots, i$}{
		$h_{mi} = \vek v_{m}^{T}\vek v_{i+1} $\\
		$\vek v_{i+1}  \leftarrow \vek v_{i+1}  - h_{mi} \vek v_{m}$
	}
	$h_{i+1,i}= \norm{\vek v_{i+1}}$\\
	$\vek v_{i+1} = \vek v_{i+1}/h_{i+1,i}$\\
	$\vek d_{i+1} = \vek v_{i+1}^{T}\vek C$\\
	$\vek m_{i+1} \leftarrow \vek d_{i+1}$\\
	Apply previous rotations to $j$th column of $\underline{\vek H}_{j}$\\
	Obtain Givens sine and cosine $s_{j}$ and $c_{j}$ and updated $\vek R_{j}$\\
	Apply new rotations to update $\widehat{\vek b}_{j} = \vek Q_{j}^{T}\vek V_{j+1}^{T}\vek r_{0}$\\
	$\begin{bmatrix}
		\vek m_{j}
		\\
		\vek m_{j+1}
	\end{bmatrix}
	\leftarrow
	\begin{bmatrix}
		c_{j}
		&
		s_{j}
		\\
		-s_{j}
		&
		c_{j}
	\end{bmatrix}
	\begin{bmatrix}
		\vek m_{j}
		\\
		\vek m_{j+1}
	\end{bmatrix}$\\
	\If{$\gamma \cdot \abs{\widehat{\vek b}_{j}(j+1)} < \nm{\vek r_{0}}\varepsilon_{tol}$ }{
		Solve $\prn{\vek I - \vek M_{j}\vek M_{j}^{T}}\vek R_{j}\vek y_{j} = \widehat{\vek b}_{j}(1:j)$\\
		Set $\vek t \leftarrow \vek V_{j}\vek y_{j}$\\
		Set $\vek s_{2} = -\vek U\prn{\vek F^{-1} \prn{\vek D_{j}^{T} \vek H_{j}\vek y_{j}}}$\\
		Set $\vek x \leftarrow \vek x_{0} + \vek s_{1} + \vek s_{2} + \vek t$; $\vek r \leftarrow \vek b - \vek A\vek x$\\
		\uIf{$\nm{\vek r}  < \nm{\vek r_{0}}\varepsilon_{tol}$}{
			Exit loop and {\tt return}
		}\uElse{
			$\gamma \leftarrow \nm{\vek r - \vek C\vek C^{T}\vek r}/\nm{\vek r}$
		}
	}
}
\ealgte

\KMS{%
\subsection{Comparison of implementations}
We compare Algorithm \ref{alg.r3gmres} to \cite[Algorithm 2]{DGH.2014} by studying their modifications to the common GMRES 
implementation upon which they are built, i.e., a Givens-rotation-based implementation as described in \cite{Saad.GMRES.1986}.
As in \cite{DGH.2014}, we consider operations occurring inside the outermost loop.
Inside of the main loop, both algorithms perform one matrix-vector product and an Arnoldi orthogonalization of each new basis vector.
At the beginning of the algorithm, they perform many of the same or comparable initialization steps.  According to the authors,
\cite[Algorithm 2]{DGH.2014} performs $2k^{3}$ operations for additional Givens rotations per iteration since that method treats the minimization 
\ref{eqn.rrrGMRES-small-ls} directly. Additionally, obtaining an update of $\vek C_{j}$ at each iteration costs $2kn$ operations, 
and obtaining $\vek F_{j}$ costs $2k^{2}n$ at each iteration.  Additionally, there are some lower-order costs.  Thus, \cite[Algorithm 2]{DGH.2014}
has a per-iteration cost above that of GMRES of roughly $2k(k^{2} + kn + n)$ operations.

The formulation of Algorithm \ref{alg.r3gmres} allows us to discard many of these per-iteration operations.  A comparable operation
which is not discarded is the progressive building of $\vek D_{j}$, which costs $nk$ operations.  The update of $\vek M_{j}$ costs $2k$ 
operations.  Thus, Algorithm \ref{alg.r3gmres} has a per-iteration cost above that of GMRES of roughly 
$k(n+2)$.  These are the per-iteration cost of both algorithms above that of GMRES is $\CO(n)$.  The main difference is that the per iteration cost
of Algorithm \ref{alg.r3gmres} above GMRES is linear in $k$ whereas it is cubic for \cite[Algorithm 2]{DGH.2014}.  
Thus we conclude that Algorithm \ref{alg.r3gmres} can accommodate a larger augmentation subspace with only linear growth in cost of additional
operations.  
}%

\section{Numerical Results}\label{section.numerics}
In this section, we demonstrate that Algorithm \ref{alg.r3gmres} \KMS{produces approximations of the same quality as those 
produced} by the version of the algorithm presented in \cite{DGH.2014} using
code from the authors.  The point here is not to compare the superiority of one version or the other, as neither code is optimized.  Rather,
as this note is laying out an alternative approach to the augmentation of unprojected Krylov subspaces, we demonstrate that our code delivers the
same performance, verifying the alternative mathematical derivation in previous sections. 
We reproduce two experiments from \cite{DGH.2014} using {\tt Regularization Tools} \cite{H.2007} with problem size $n=256$.
The noise vectors are generated from the normal distribution using {\tt randn()}.  For the experiments, we report the level of the noise relative to the
size of the right-hand side, i.e., a relative noise level of $10^{-3}$ means that the $2$-norm of the vector perturbing the right-hand-side $\vek b_{true}$
is $10^{-3}\nm{\vek b_{true}}$.
All experiments are performed in Matlab R2020a and we have established a repository \cite{Soodhalter:2021:1} in which our code for
Algorithm \ref{alg.r3gmres} is contained. 

\subsection{Experiment: {\tt deriv2()} test}
\begin{figure}
\caption{Experiment: {\tt deriv2()} test\label{fig.deriv2}}
\includegraphics[scale=.2]{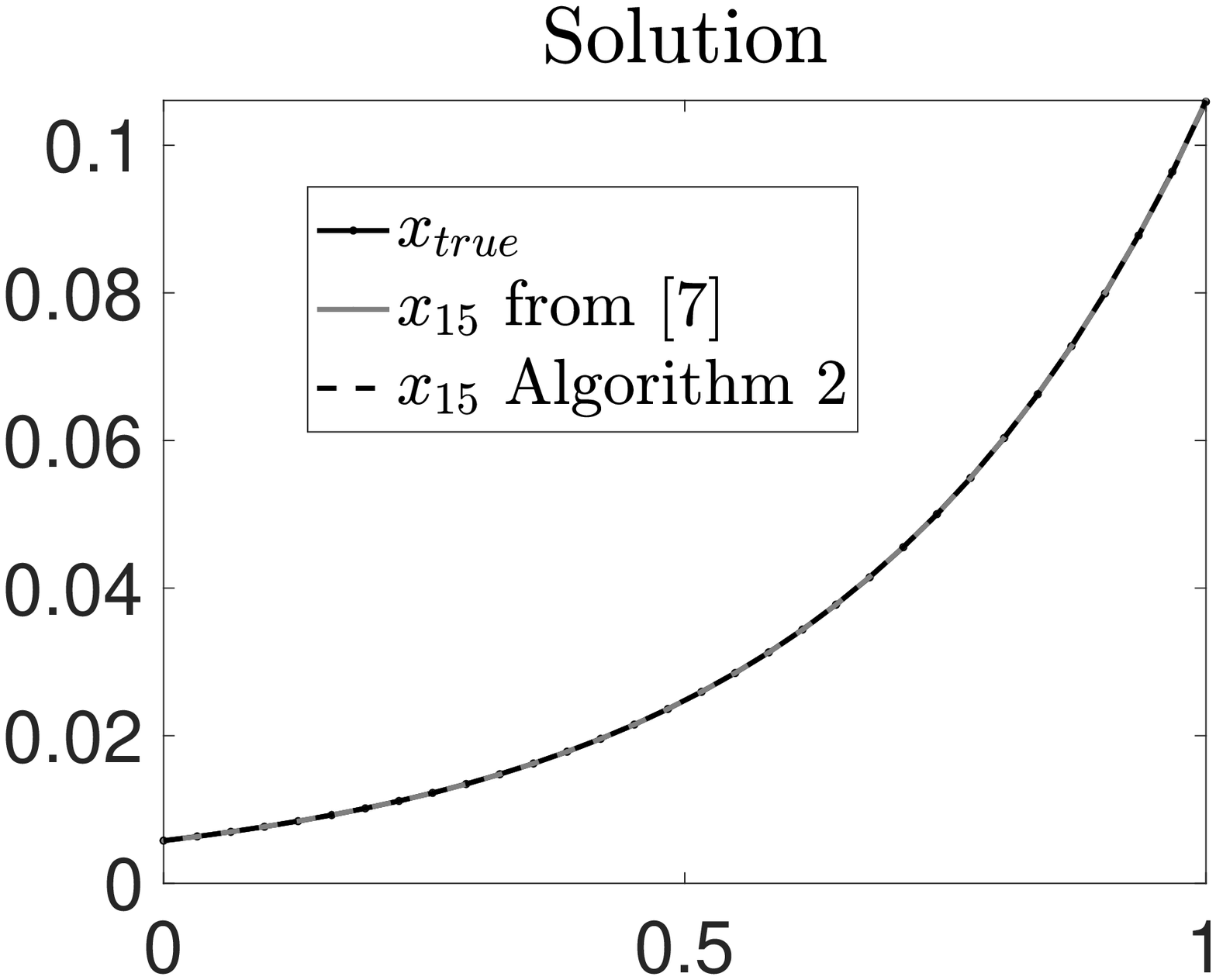}
\includegraphics[scale=.2]{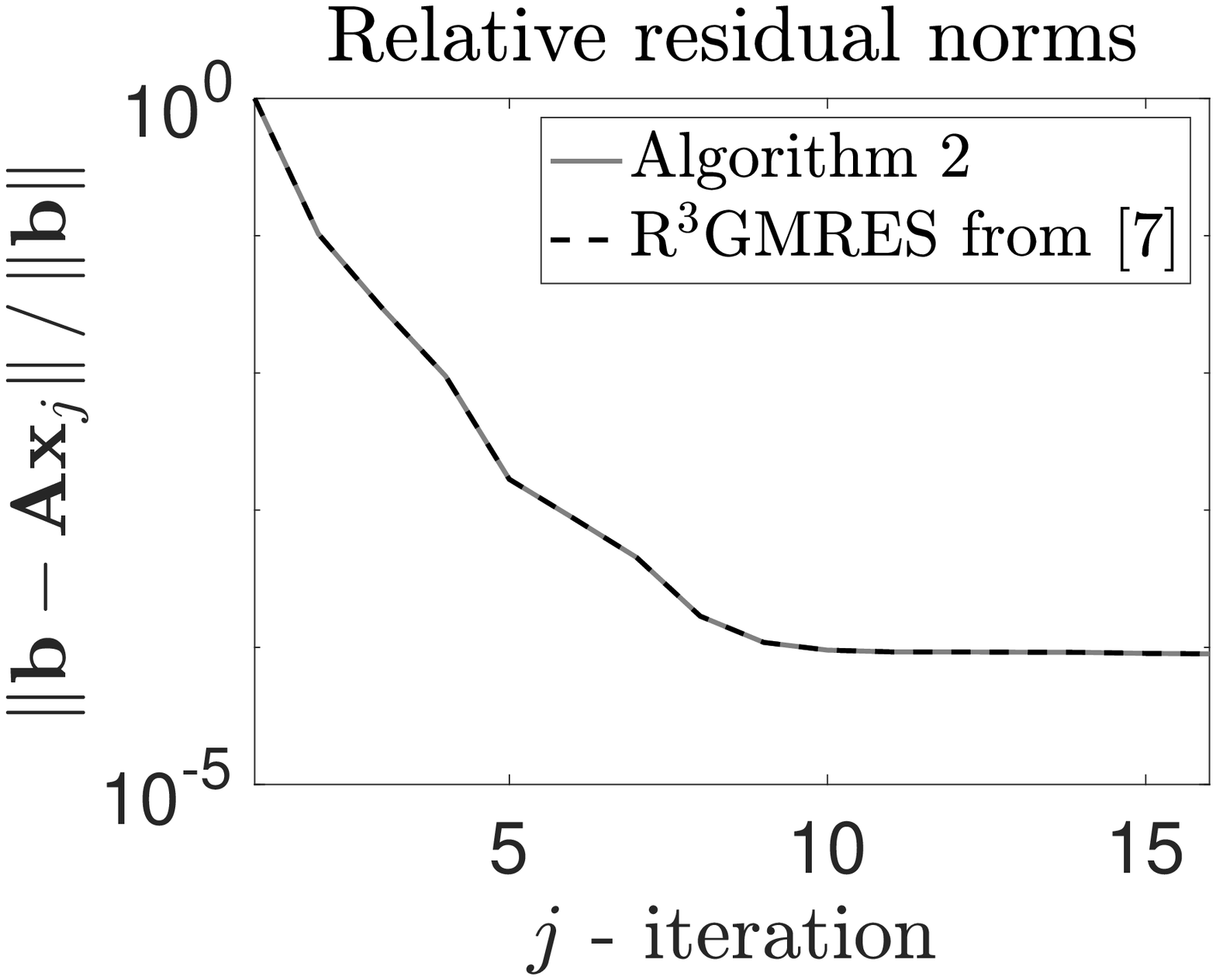}
\includegraphics[scale=.2]{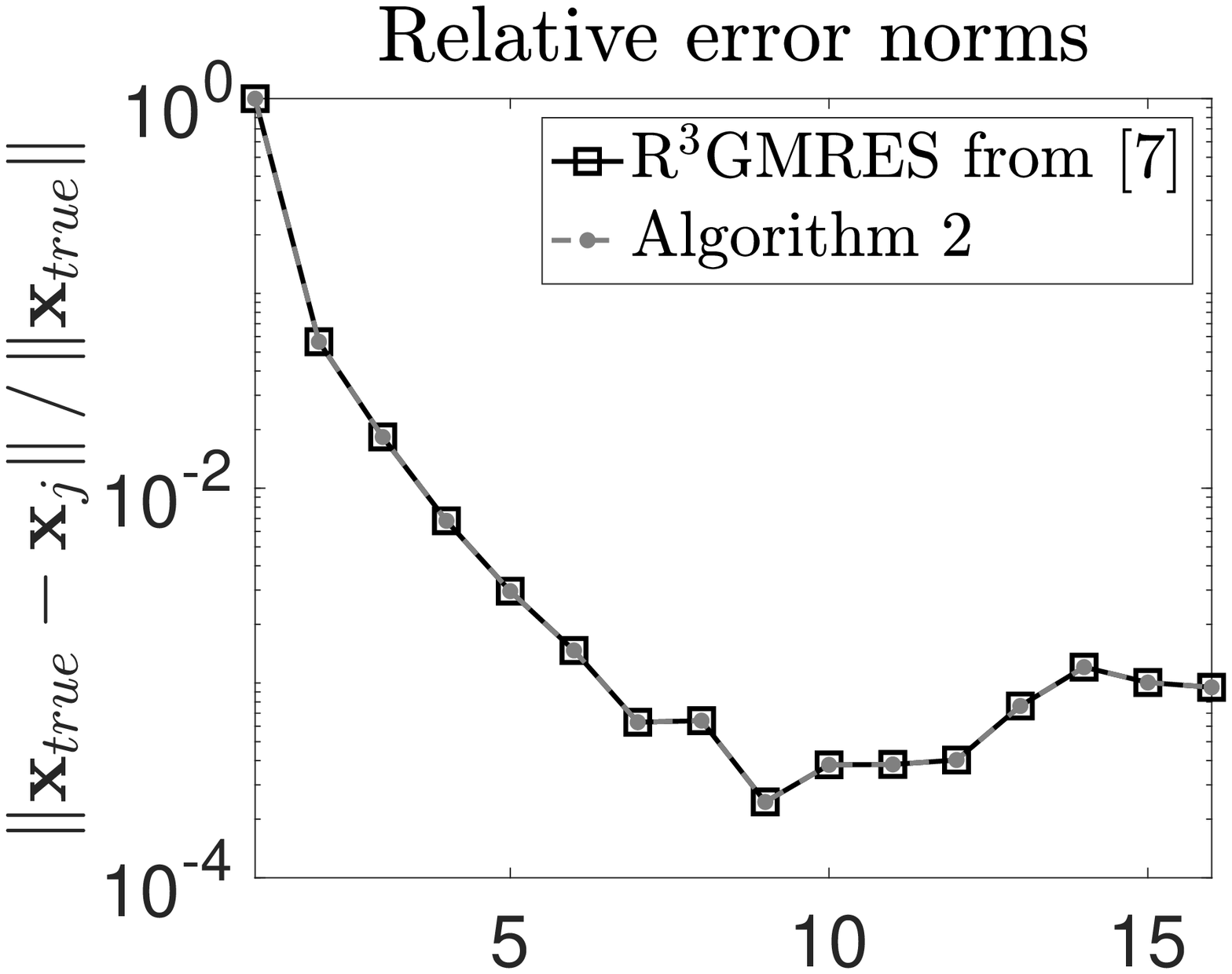}
\end{figure}

This reproduces the experiment in \cite[Section 4.2]{DGH.2014} wherein augmentation is used to help encode known boundary conditions approximately so
that the iteration focuses mostly on reconstructing the solution on the interior of the domain.  The matrix is generated by the {\tt deriv2()} function which
produces a discretization of the Fredholm integral operator whose kernel is the Green's function of the second derivative operator.
The relative noise level is $10^{-5}$.
Following \cite[Section 4.2]{DGH.2014}, we set
$\CU=\Span\curl{\bbmat 1 & 1 & \cdots & 1 \ebmat^{T}, \bbmat 1 & 2 & \cdots & n \ebmat^{T}}$.   
Results shown in \Cref{fig.deriv2} demonstrate that the performance of
the two implementations is virtually indistinguishable. 

\subsection{Experiment: {\tt gravity()} test -- correctly localized discontinuity}
\begin{figure}
\caption{Experiment: {\tt gravity()} test -- correctly localized discontinuity\label{fig.gravity}}
\includegraphics[scale=.2]{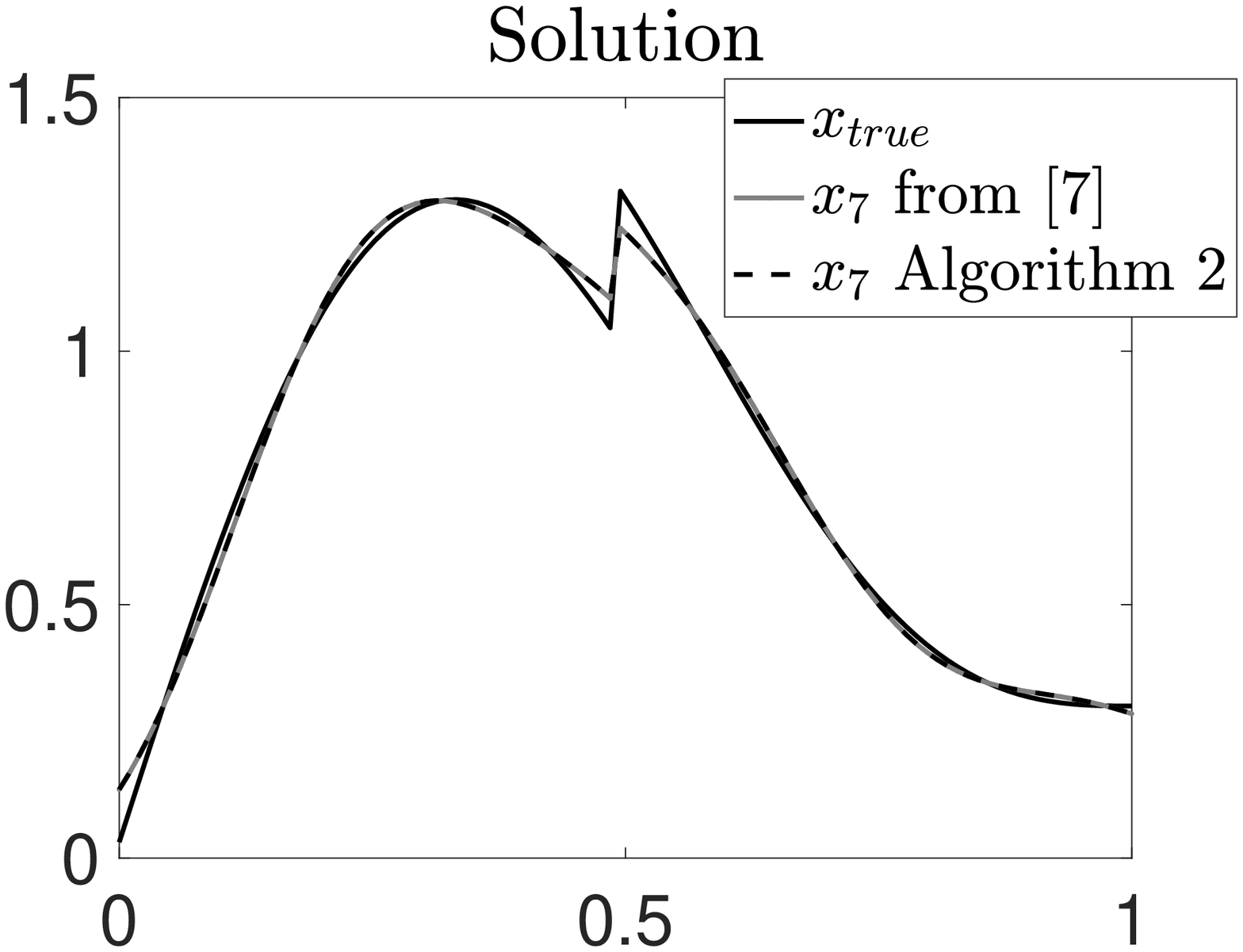}
\includegraphics[scale=.2]{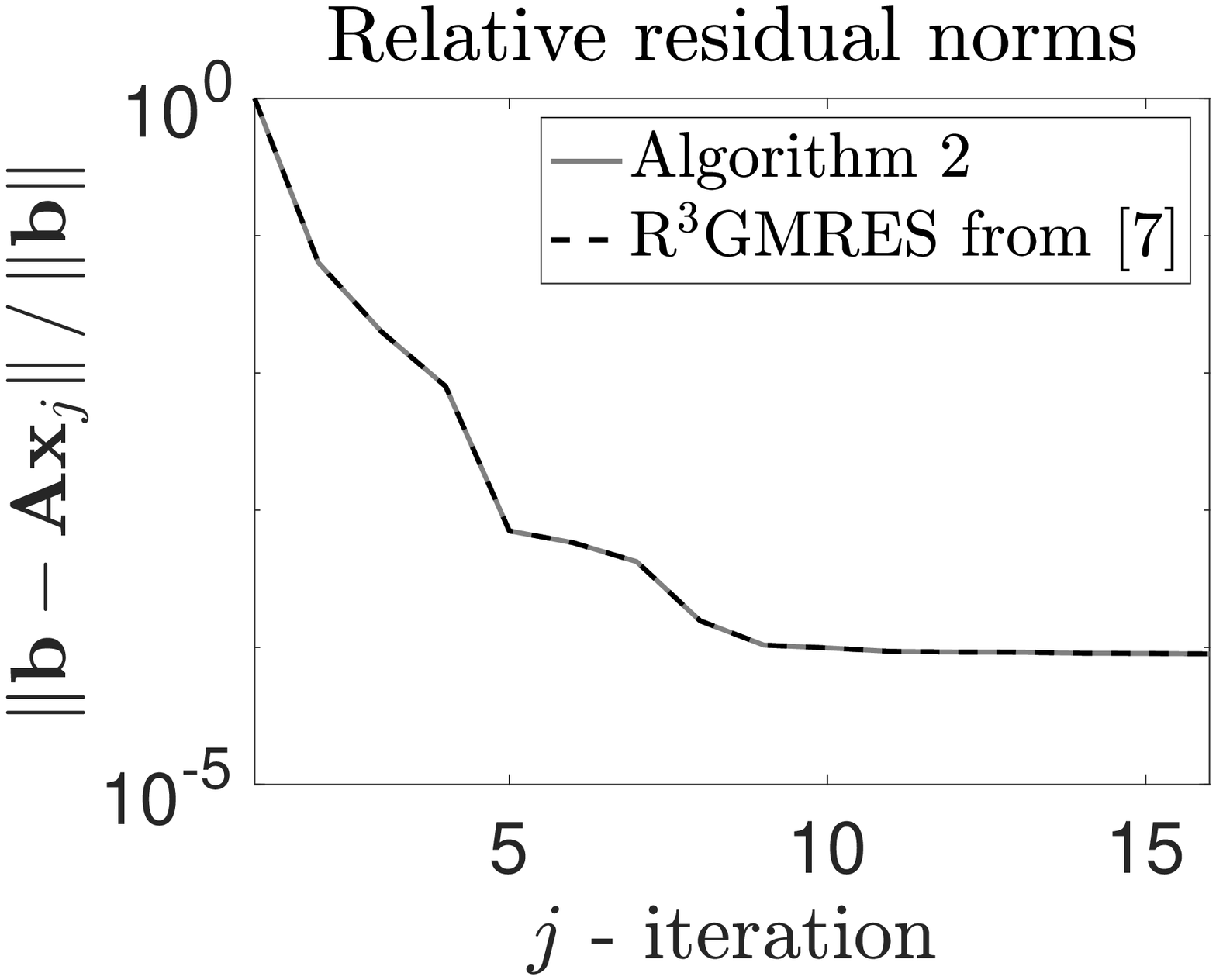}
\includegraphics[scale=.2]{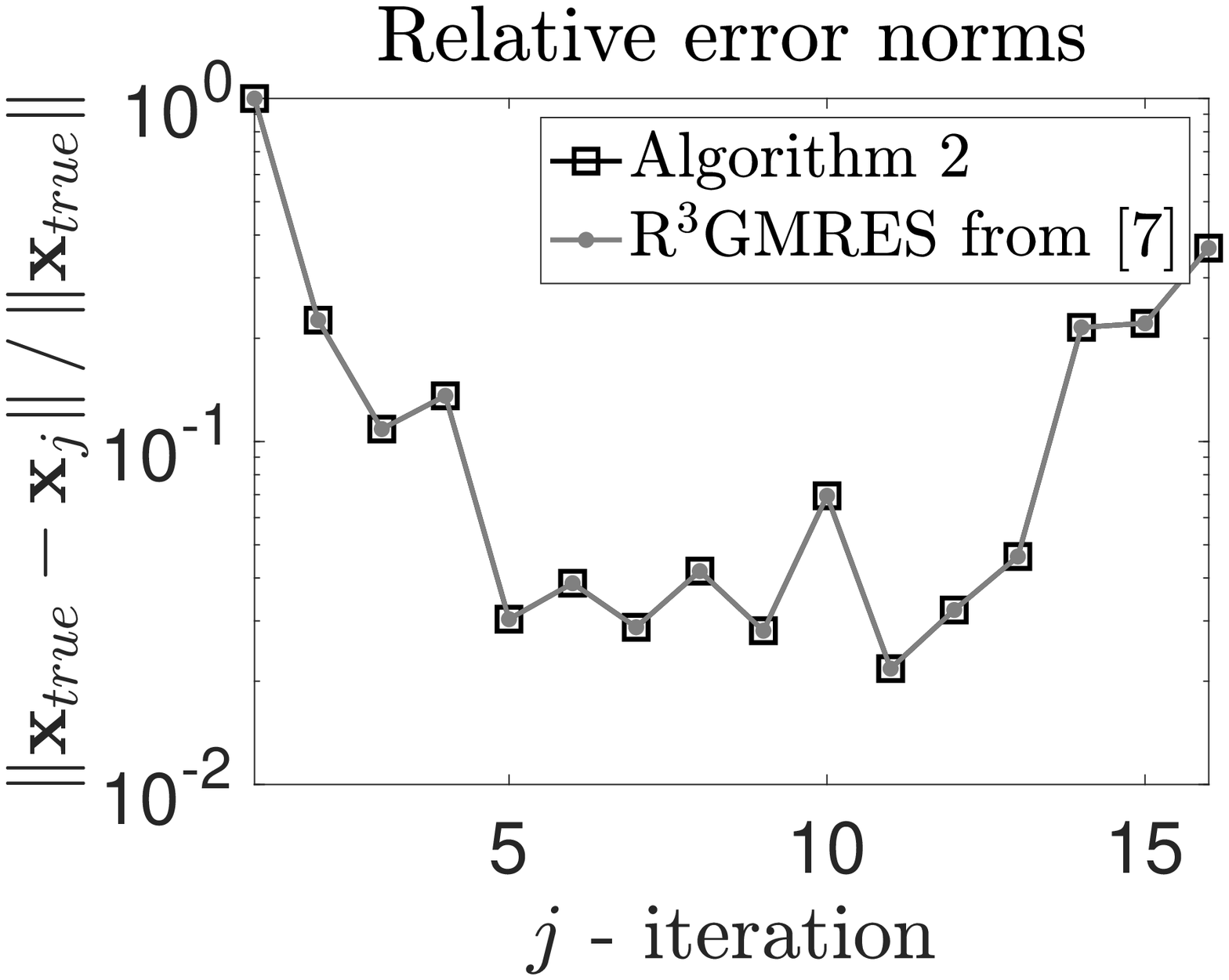}
\end{figure}

We generate the matrix $\vek A$ for this example using the {\tt gravity()} function, which generates a discretization of a 
Fredholm integral operator of the first kind modeling a one-dimensional gravity surveying problem application posed on the interval $[0,1]$.  
Relative noise level is $10^{-4}$.
We take the true solution produced by the function and introduce a discontinuity at $t=\dfrac{1}{2}$, as in \cite[Section 4.3]{DGH.2014}. For this 
experiment, we assume we know the location of the discontinuity and set $\CU=\Span\curl{\bbmat 0 & 0 & \cdots & 0 & 1 & \cdots & 1 \ebmat^{T}}$
to correctly encode this discontinuity.  In \Cref{fig.gravity}, we see that the two implementations perform identically.

\subsection{Experiment: {\tt gravity()} test -- incorrectly localized discontinuity}
\begin{figure}
\caption{Experiment: {\tt gravity()} test -- incorrectly localized discontinuity\label{fig.gravity-baddiscont}}
\includegraphics[scale=.2]{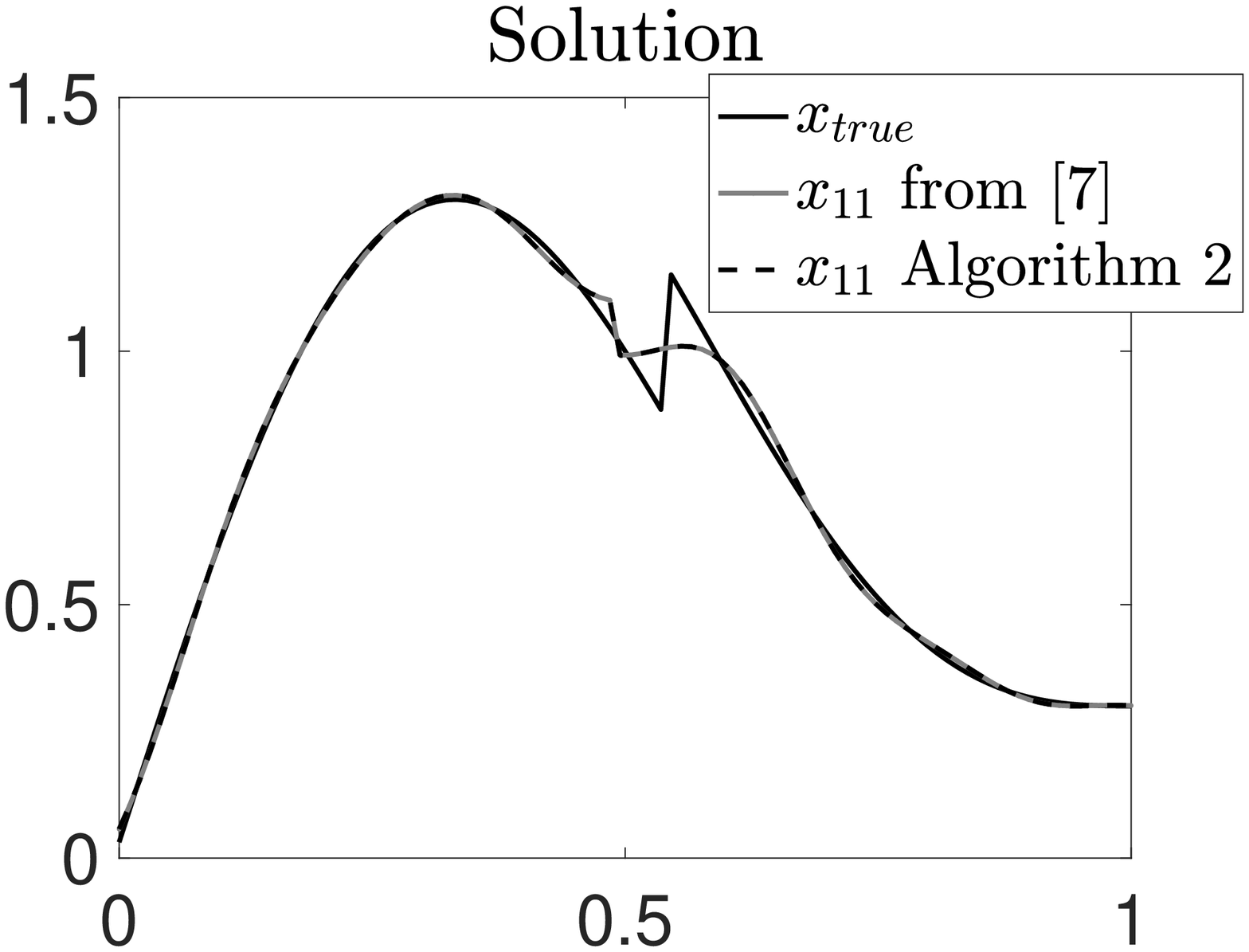}
\includegraphics[scale=.2]{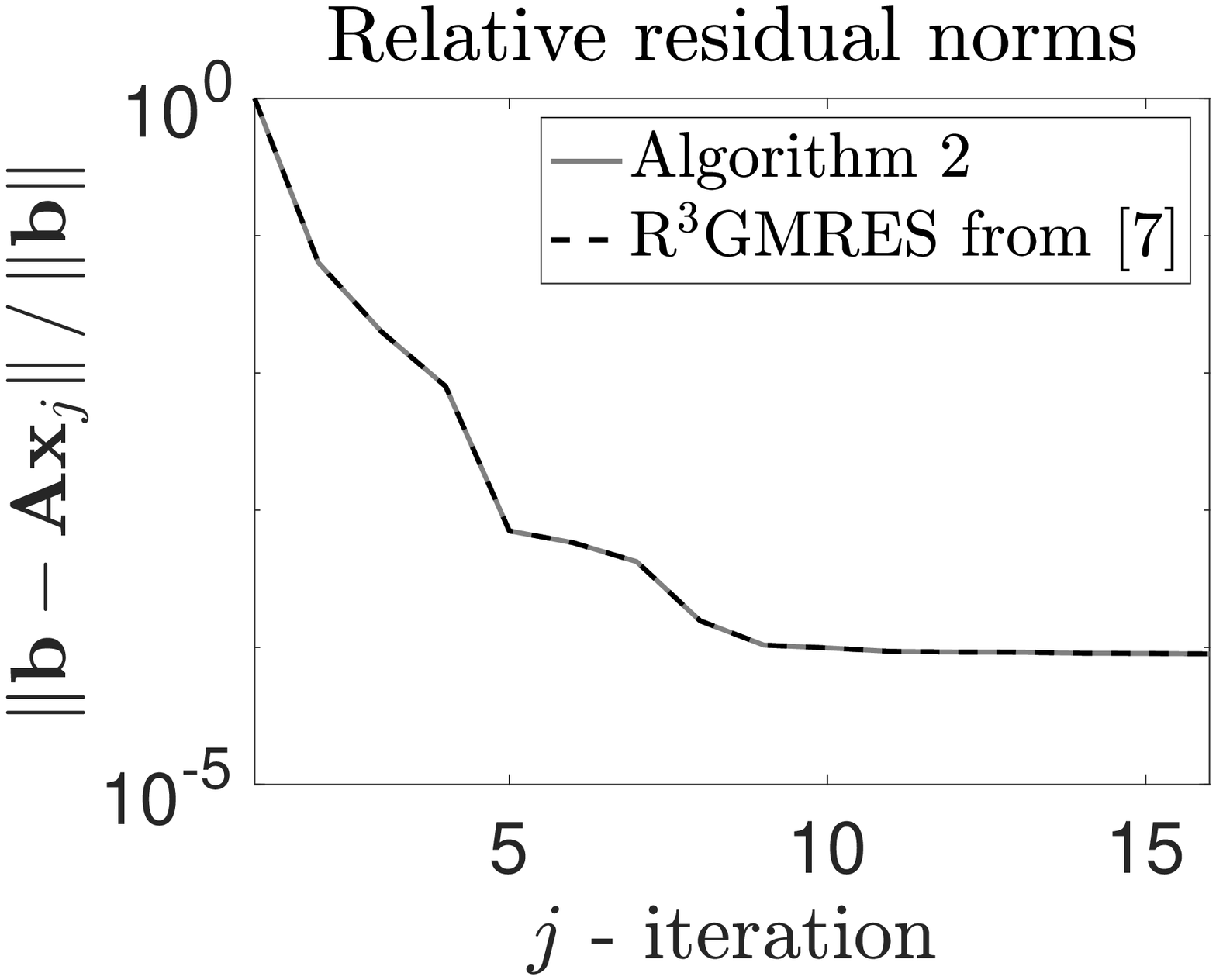}
\includegraphics[scale=.2]{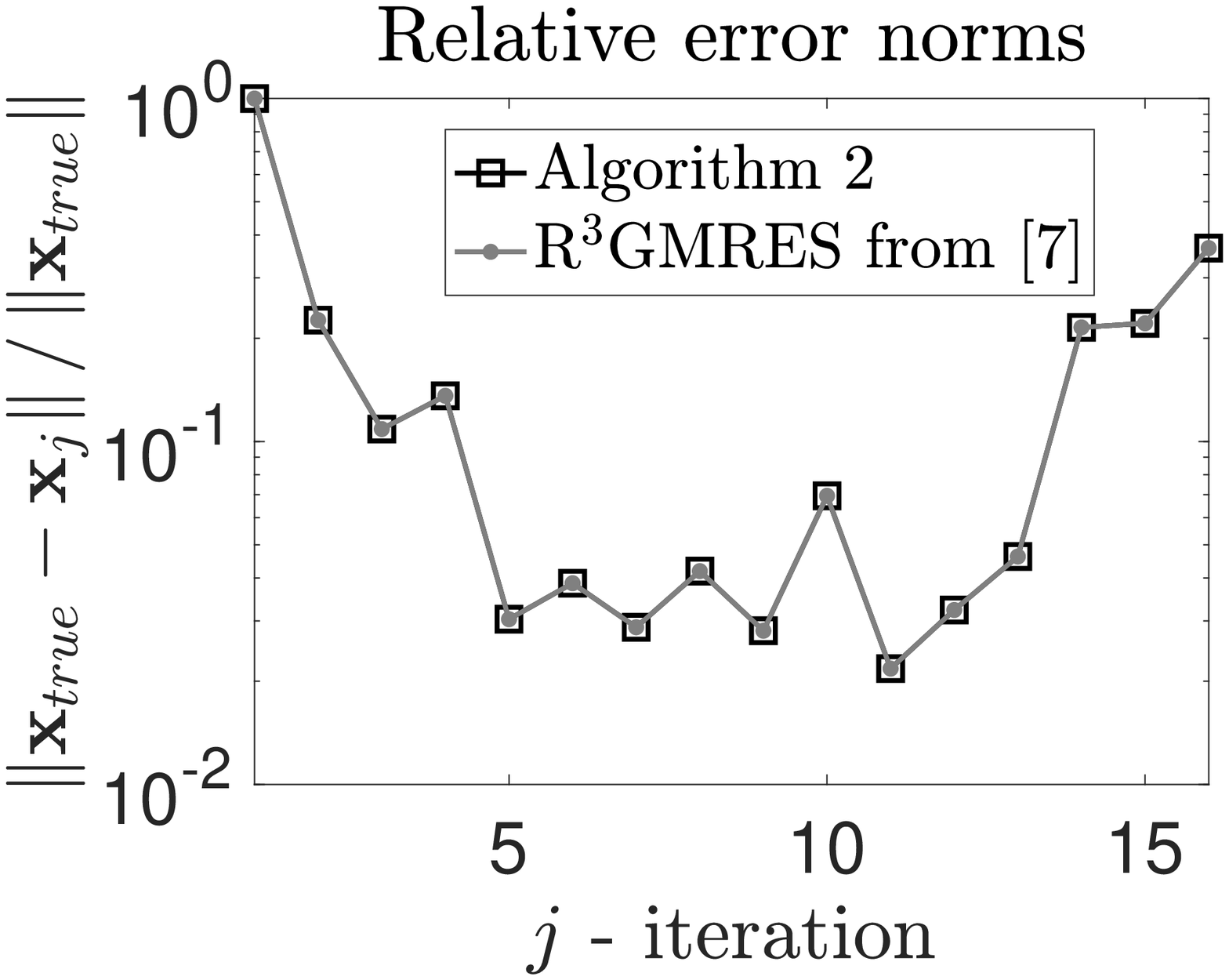}
\end{figure}

\begin{figure}
\centering
\caption{\label{fig.gravity-baddiscont-gcrodr}  {\tt gravity()} test -- incorrectly localized discontinuity with projected Krylov subspace}
\includegraphics[scale=.3]{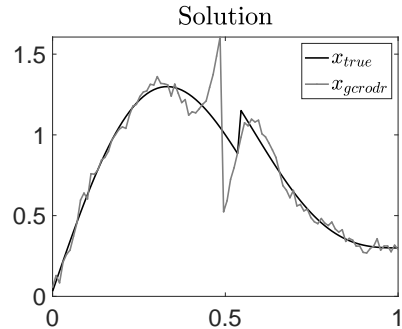}
\end{figure}

We construct the same problem as in the previous experiment, but we move the discontinuity to a $t>\dfrac{1}{2}$.  However, we encode the discontinuity
incorrectly using the same $\CU$ as in the previous experiment. In \Cref{fig.gravity-baddiscont}, we observe that both implementations again perform identically.
Furthermore, one sees that \KMS{the} minimization process reduces the influence of the falsely-placed discontinuity encoded by $\CU$ while trying to fit the true discontinuity.
This has been noted in \cite[Section 4.3]{DGH.2014} as a possible advantage in augmenting an unprojected Krylov subspace for solving an ill-posed problem, as the 
incorrectly-chosen $\CU$ does not influence which Krylov subspace is built.
We contrast this with the best reconstruction produced by an augmented iterative solver using a projected Krylov subspace, $\CK\prn{ \prn{\vek I-\bPhi}\vek A, \prn{\vek I-\bPhi}\vek r_{0}}$, using a GCRO-type code.  In \Cref{fig.gravity-baddiscont-gcrodr}, we see that the method at its best still emphasizes the incorrectly
localized discontinuity.

\section{Discussion}\label{section.conclusions}
The main goal in this note is to demonstrate that augmented \KMS{unprojected} Krylov subspace methods fit into the same framework from \cite{dSKS.2020} enabling a
simpler implementation in the style of a GCRO-DR type method.  This leads us to observe that the R$^{3}$GMRES method is closely related
to the augmentation strategy from \cite{Saad.Deflated-Aug-Krylov.1997}.  
With that perspective, we show one can actually approximate the solution to a projected subproblem and project
the approximation to obtain the part from the augmented subspace.  The benefit when applying this to the R$^{3}$GMRES method is we no longer need
to progressively maintain an orthonormal basis to the full sum subspace $\vek A\prn{\CU + \CV_{j}}$.  

The numerical experiments we showed follow from
what was done in \cite{DGH.2014}, focusing on instances wherein one wants to enforce that the solution has an a priori known structure but accommodate the
possibility that this knowledge is flawed.  We contrasted this with the performance of a GCRO-type method to show how for an ill-posed problem, an augmented
method with a projected Krylov subspace can over-emphasize the bad knowledge to an extent that it cannot recover due to the ill-posedness of the problem.

However, it should be noted that GCRO-based augmentation/recycling methods still exhibit superior performance when it comes to the acceleration
of convergence for
complicated, large-scale problems.  
Rather, this work highlights that it can be important to distinguish between
``trustworthy'' and ``untrustworthy'' information when using augmentation \KMS{methods, particularly for ill-posed problems}.
A future path to explore would be to consider mixing the two strategies more generally for situations
when one has both trustworthy and untrustworthy/corrupted information one wishes to use without it corrupting the behavior of the solver.

\section*{Acknowledgments}
The author thanks Per Christian Hansen for an interesting discussion about this topic back when we could go to conferences in person and for sending the author
an implementation of R$^{3}$GMRES to validate against.  \KMS{The author also thanks the two anonymous referees for their helpful comments and
suggested edits to tighten up the exposition of the manuscript.}

\bibliographystyle{siamplain}
\bibliography{master,local}

\begin{thebibliography}{10}

\bibitem{BR-2.2007}
{\sc J.~Baglama and L.~Reichel}, {\em Augmented {GMRES}-type methods},
  Numerical Linear Algebra with Applications, 14 (2007), pp.~337--350,
  \url{https://doi.org/10.1002/nla.518},
  \url{http://dx.doi.org/10.1002/nla.518}.

\bibitem{BR.2007}
{\sc J.~Baglama and L.~Reichel}, {\em Decomposition methods for large linear
  discrete ill-posed problems}, Journal of Computational and Applied
  Mathematics, 198 (2007), pp.~332--343,
  \url{https://doi.org/10.1016/j.cam.2005.09.025},
  \url{http://dx.doi.org/10.1016/j.cam.2005.09.025}.

\bibitem{ChungGazzola:2021:1}
{\sc J.~Chung and S.~Gazzola}, {\em Computational methods for large-scale
  inverse problems: a survey on hybrid projection methods},
  \url{https://arxiv.org/abs/2105.07221}.

\bibitem{deSturler.GCRO.1996}
{\sc E.~{de Sturler}}, {\em Nested {K}rylov methods based on {GCR}}, Journal of
  Computational and Applied Mathematics, 67 (1996), pp.~15--41,
  \url{https://doi.org/10.1016/0377-0427(94)00123-5},
  \url{http://dx.doi.org/10.1016/0377-0427(94)00123-5}.

\bibitem{deSturler.GCROT.1999}
{\sc E.~{de Sturler}}, {\em Truncation strategies for optimal {K}rylov subspace
  methods}, SIAM Journal on Numerical Analysis, 36 (1999), pp.~864--889,
  \url{https://doi.org/10.1137/S0036142997315950},
  \url{http://dx.doi.org/10.1137/S0036142997315950}.

\bibitem{EdS-personal.2020}
{\sc E.~de~Sturler}.
\newblock private communication, 2020.

\bibitem{DGH.2014}
{\sc Y.~Dong, H.~Garde, and P.~C. Hansen}, {\em R{${}^3$}{GMRES}: including
  prior information in {GMRES}-type methods for discrete inverse problems},
  Electron. Trans. Numer. Anal., 42 (2014), pp.~136--146.

\bibitem{DMR.2014}
{\sc L.~Dykes, F.~Marcell\'{a}n, and L.~Reichel}, {\em The structure of
  iterative methods for symmetric linear discrete ill-posed problems}, BIT, 54
  (2014), pp.~129--145, \url{https://doi.org/10.1007/s10543-014-0476-2},
  \url{https://doi-org.libproxy.temple.edu/10.1007/s10543-014-0476-2}.

\bibitem{ErhelGuyomarch:2000:1}
{\sc J.~Erhel and F.~Guyomarc'h}, {\em An augmented conjugate gradient method
  for solving consecutive symmetric positive definite linear systems}, 21,
  pp.~1279--1299, \url{https://doi.org/10.1137/s0895479897330194},
  \url{https://doi.org/10.1137%2Fs0895479897330194}.

\bibitem{Gaul.2014-phd}
{\sc A.~Gaul}, {\em Recycling {K}rylov subspace methods for sequences of linear
  systems: Analysis and applications}, PhD thesis, Technischen Universit\"at
  Berlin, 2014.

\bibitem{GGL.2013}
{\sc A.~Gaul, M.~H. Gutknecht, J.~Liesen, and R.~Nabben}, {\em A framework for
  deflated and augmented {K}rylov subspace methods}, SIAM Journal on Matrix
  Analysis and Applications, 34 (2013), pp.~495--518,
  \url{https://doi.org/10.1137/110820713},
  \url{http://dx.doi.org/10.1137/110820713}.

\bibitem{GazzolaLandman:2020:1}
{\sc S.~Gazzola and M.~S. Landman}, {\em Krylov methods for inverse problems:
  Surveying classical, and introducing new, algorithmic approaches}, 43,
  \url{https://doi.org/10.1002/gamm.202000017},
  \url{https://doi.org/10.1002%2Fgamm.202000017}.

\bibitem{Gutknecht:2012:1}
{\sc M.~H. Gutknecht}, {\em Spectral deflation in {K}rylov solvers: a theory of
  coordinate space based methods}, Electron. Trans. Numer. Anal., 39 (2012),
  pp.~156--185.

\bibitem{Gutknecht.AugBiCG.2014}
{\sc M.~H. Gutknecht}, {\em Deflated and augmented {K}rylov subspace methods: a
  framework for deflated {B}i{CG} and related solvers}, SIAM J. Matrix Anal.
  Appl., 35 (2014), pp.~1444--1466, \url{https://doi.org/10.1137/130923087},
  \url{https://doi-org.libproxy.temple.edu/10.1137/130923087}.

\bibitem{H1995-book}
{\sc M.~Hanke}, {\em Conjugate gradient type methods for ill-posed problems},
  vol.~327 of Pitman Research Notes in Mathematics Series, Longman Scientific
  \& Technical, Harlow, 1995.

\bibitem{H.2007}
{\sc P.~C. Hansen}, {\em Regularization tools version 4.0 for matlab 7.3},
  Numerical Algorithms, 46 (2007), pp.~189--194.

\bibitem{dSC.2019}
{\sc J.~Jiang, J.~Chung, and E.~de~Sturler}, {\em Hybrid projection methods
  with recycling for inverse problems}, pp.~S146--S172,
  \url{https://doi.org/10.1137/20m1349515},
  \url{https://doi.org/10.1137%2F20m1349515}.

\bibitem{Morgan.GMRESDR.2002}
{\sc R.~B. Morgan}, {\em G{MRES} with deflated restarting}, SIAM Journal on
  Scientific Computing, 24 (2002), pp.~20--37,
  \url{https://doi.org/10.1137/S1064827599364659},
  \url{http://dx.doi.org/10.1137/S1064827599364659}.

\bibitem{NRS-2.2012}
{\sc A.~Neuman, L.~Reichel, and H.~Sadok}, {\em Algorithms for range restricted
  iterative methods for linear discrete ill-posed problems}, Numer. Algorithms,
  59 (2012), pp.~325--331, \url{https://doi.org/10.1007/s11075-011-9491-4},
  \url{http://dx.doi.org/10.1007/s11075-011-9491-4}.

\bibitem{NRS.2012}
{\sc A.~Neuman, L.~Reichel, and H.~Sadok}, {\em Implementations of range
  restricted iterative methods for linear discrete ill-posed problems}, Linear
  Algebra Appl., 436 (2012), pp.~3974--3990,
  \url{https://doi.org/10.1016/j.laa.2010.08.033},
  \url{http://dx.doi.org/10.1016/j.laa.2010.08.033}.

\bibitem{NeumanReichelSadok:2012}
{\sc A.~Neuman, L.~Reichel, and H.~Sadok}, {\em Implementations of range
  restricted iterative methods for linear discrete ill-posed problems}, Linear
  Algebra and its Applications, 436 (2012), pp.~3974--3990,
  \url{https://doi.org/10.1016/j.laa.2010.08.033},
  \url{https://doi.org/10.1016%2Fj.laa.2010.08.033}.

\bibitem{Parks.deSturler.GCRODR.2005}
{\sc M.~L. Parks, E.~de~Sturler, G.~Mackey, D.~D. Johnson, and S.~Maiti}, {\em
  Recycling {K}rylov subspaces for sequences of linear systems}, SIAM Journal
  on Scientific Computing, 28 (2006), pp.~1651--1674,
  \url{https://doi.org/10.1137/040607277},
  \url{http://dx.doi.ofrrg/10.1137/040607277}.

\bibitem{Parks.Soodhalter.Szyld.16}
{\sc M.~L. Parks, K.~M. Soodhalter, and D.~B. Szyld}, {\em A block recycled
  gmres method with investigations into aspects of solver performance},
  \url{https://arxiv.org/abs/1604.01713}.

\bibitem{HRS.Augmented-regularization.2020}
{\sc R.~Ramlau, K.~M. Soodhalter, and V.~Hutterer}, {\em Subspace
  recycling-based regularization methods},
  \url{https://arxiv.org/abs/2011.05473}.

\bibitem{RamlauStadler:2020:1}
{\sc R.~Ramlau and B.~Stadler}, {\em An augmented wavelet reconstructor for
  atmospheric tomography}, \url{https://arxiv.org/abs/2011.06842}.

\bibitem{RY.2005}
{\sc L.~Reichel and Q.~Ye}, {\em Breakdown-free {GMRES} for singular systems},
  SIAM J. Matrix Anal. Appl., 26 (2005), pp.~1001--1021 (electronic),
  \url{https://doi.org/10.1137/S0895479803437803},
  \url{http://dx.doi.org/10.1137/S0895479803437803}.

\bibitem{Saad.Deflated-Aug-Krylov.1997}
{\sc Y.~Saad}, {\em Analysis of augmented {K}rylov subspace methods}, SIAM
  Journal on Matrix Analysis and Applications, 18 (1997), pp.~435--449,
  \url{https://doi.org/10.1137/S0895479895294289},
  \url{http://dx.doi.org/10.1137/S0895479895294289}.

\bibitem{Saad.GMRES.1986}
{\sc Y.~Saad and M.~H. Schultz}, {\em {GMRES}: A generalized minimal residual
  algorithm for solving nonsymmetric linear systems}, SIAM Journal on
  Scientific and Statistical Computing, 7 (1986), pp.~856--869.

\bibitem{Soodhalter:2021:1}
{\sc K.~M. Soodhalter}, {\em kirkmsoodhalter/r3gmres-simplified: R3gmres
  simplified implementation}, \url{https://doi.org/10.5281/ZENODO.4975990},
  \url{https://zenodo.org/record/4975990}.

\bibitem{dSKS.2020}
{\sc K.~M. Soodhalter, E.~d. Sturler, and M.~E. Kilmer}, {\em A survey of
  subspace recycling iterative methods}, GAMM-Mitt., 43 (2020),
  pp.~e202000016--28.

\end{thebibliography}
\end{document}